\documentclass[11pt]{article}
\usepackage{graphicx}
\usepackage{subfigure}
\usepackage{color}
\usepackage{amsmath,amssymb}%,latexsym,amsthm}
\def\la{\Big\langle}
\def\ra{\Big\rangle}

\def\forall{\hbox{for all}~~}
\def\bfe{{\bf e}}
\def\bfk{{\bf k}}
\def\L{{\bf L}}

\def\E{{\cal E}}

\def\bfv{{\bf v}}

\def\bfb{{\bf b}}
\def\bfc{{\bf c}}
\def\bfh{{\bf h}}
\def\bfz{{\bf z}}
\def\bfw{{\bf w}}

\def\bfp{{\bf p}}
\def\ds{\displaystyle}
\def\cL{{\cal L}}

\def\bfn{{\bf n}}
\def\ve{\varepsilon}

\def\D{{\cal D}}

\def\V{{\cal V}}

\def\A{{\cal A}}

\def\R{{\mathbb R}}

\def\implies{\Longrightarrow}
\def\vp{\varphi}

\def\P{{\cal P}}

\def\v{\vskip 1em}
\def\O{{\cal O}}

\def\C{{\cal C}}
\def\weakto{\rightharpoonup}
\def\ov{\overline}
\def\Tilde{\widetilde}

\def\bel{\begin{equation}\label}
\def\eeq{\end{equation}}
\def\sqr#1#2{\vbox{\hrule height .#2pt
\hbox{\vrule width .#2pt height #1pt \kern #1pt
\vrule width .#2pt}\hrule height .#2pt }}
\def\square{\sqr74}
\def\endproof{\hphantom{MM}\hfill\llap{$\square$}\goodbreak}
\def\bega{\begin{array}}
\def\enda{\end{array}}
\def\begi{\begin{itemize}}
\def\endi{\end{itemize}}

\begin{document}
\title{\bf  Growth Models for Tree Stems and Vines}

\author{Alberto Bressan, Michele Palladino, and Wen Shen\\ \, \\
Department of Mathematics, Penn State University.\\
University Park, PA~16802, USA.\\
\\
e-mails:~axb62@psu.edu,
~mup26@psu.edu,  
~wxs27@psu.edu}
\maketitle
\begin{abstract} The paper introduces a PDE model for the growth 
of a tree stem or a vine.   The equations
describe the elongation due to cell growth, and the response to gravity 
and to external obstacles. An additional term accounts 
for the tendency of a vine to curl around branches of other plants.

When obstacles are present, 
the model takes  the form of a differential inclusion 
with unilateral constraints.
At each time $t$, a cone of admissible reactions is determined by
the minimization of an elastic deformation energy.
The main theorem shows that local solutions exist and can be prolonged 
globally in time, except when
a specific ``breakdown configuration" is reached.
Approximate solutions are constructed by an
operator-splitting technique.   
Some numerical simulations
are provided at the end of the paper.
\end{abstract}

\section{Introduction}
\label{s:0}
\setcounter{equation}{0}

We consider a simple mathematical model describing how 
the stem of a plant grows, and how it reacts to external 
constraints, 
such as branches of other plants.
At each time $t$ the stem is described by a curve $\gamma(t,\cdot)$
in 3-dimensional space.  The model takes into account
the linear elongation due to cell growth and the upward bending
as a response to gravity.   In the case of vines, 
an additional term accounts 
for the tendency to curl around branches of other plants.

{}From a theoretical perspective, the main challenge comes from the presence
of external obstacles, %such as rocks or branches of other plants,
resulting in a number of unilateral constraints.
Ultimately, this yields a differential inclusion on a closed subset 
of  $H^2([0,T];\,\R^3)$. 
We remark that most of the literature on differential inclusions
with constraints
is concerned with the case where the cone of admissible 
reactions produced by 
the (possibly moving) obstacle is perpendicular to its boundary
\protect\cite{CG, CMM, M, RS}. In Moreau's ``sweeping process", this assumption 
plays an essential role in the proof 
of existence  and continuous dependence of solutions.
In our model, at a time when part of the stem touches the obstacle,
the evolution is governed by
the minimization of an instantaneous elastic deformation 
energy, subject to the external constraints. As a consequence, 
the cone of
admissible velocities determined by the obstacle's reaction
can be very different from the normal cone. 
In certain  ``breakdown configurations", as 
shown in Fig.~\ref{f:sg98}
this cone of admissible reactions actually happens to be tangent. 
 
Our main result, Theorem 1 in Section~\ref{s:2}, 
establishes the local existence of solutions to 
the growth model with obstacles. These solutions can be 
extended globally in time, provided that 
a specific ``breakdown configuration" is never reached. 
As already mentioned, since the cone of admissible reactions
is not a normal cone, the uniqueness and continuous dependence
of solutions is a difficult problem that requires a substantially different 
approach from \protect\cite{CG, CMM, M, RS}. A detailed analysis
of this issue will appear in the forthcoming paper \protect\cite{BP}.

The remainder of this paper is organized as follows.
Section~\ref{s:1} introduces the basic model and derives an
evolution equation satisfied by the growing curve.
If obstacles are present, this takes the form of a differential
inclusion in the space $ H^2([0,T]\,;~\R^3)$.
This is supplemented by unilateral constraints, requiring 
that at all times the curve $\gamma(t,\cdot)$ remains outside a given set.
In Section~\ref{s:2} we give a definition of solution and  state the main
existence theorem.  Namely, solutions exist locally in time
and can be prolonged up to the first time when a ``breakdown configuration" is reached.   A precise 
definition of these ``bad" configurations 
is given at (\ref{bad1})-(\ref{bad2}) and illustrated in Fig.~\ref{f:sg55}.
In essence, this happens when the tip of the stem touches the obstacle 
perpendicularly, and all the  portions of the stem that do not touch 
the obstacle are straight segments.  

The existence of solutions is proved in Sections \ref{s:44} and \ref{s:4},
constructing a sequence of approximations 
by an operator-splitting technique.   Each time step involves:
\begi
\item a regular evolution operator, modeling the linear growth and the 
bending in response to gravity (possibly including also 
the curling of vines around branches of other plants),
\item a singular operator, accounting for the obstacle reaction.
\endi
Much of the analytical work is carried out in Section~\ref{s:44}, where 
we introduce a ``push-out" operator and 
derive some key a priori estimates. 
Section~\ref{s:4} completes the proof of the main theorem.
This is based on a compactness 
argument, which yields a convergent subsequence of approximate solutions.

In Section~\ref{s:6} we briefly describe how our results
can be extended to more general models, including the case where
the elastic energies associated with twisting and bending of the stem
come with different coefficients.
Finally, Section~\ref{s:7}
presents some numerical simulations, in the case of one or two obstacles,
in two space dimensions.
The code used for these simulations can be downloaded at \protect\cite{WSweb}.
\v
\section{The basic model}
\label{s:1}
\setcounter{equation}{0}

We assume that new cells are generated at the tip of the stem, then 
they grow in size.
At time $t\geq 0$, 
the length of the cells born during the time interval
$[s, \, s+ds]$ is measured by
\bel{dL}d\ell~=~(1-e^{- \alpha(t-s)})\, ds\,,\eeq
for some constant $ \alpha>0$.
The total length of the stem  is thus
\bel{L}L(t)~
=~\int_0^t (1-e^{- \alpha(t-s)})\, ds~=~t-{1-e^{- \alpha t}\over \alpha}
\,.\eeq
%Throughout the following, for simplicity we consider the limit 
%$\alpha\to +\infty$, so that $d\ell=ds$.

At any  given time $t$, the stem will be described by a curve 
$s\mapsto P(t,s)$ in 3-dimensional space. 
For $s\in [0,t]$, the point $P(t,s)$ describes the position at time $t$
of the cell born at time $s$. 
In addition,
 we denote by $\bfk(t,s)$ the unit tangent vector to the
stem at the point $P(t,s)$, so that
\bel{k}\bfk(t,s)~=~{P_s(t,s)\over |P_s(t,s)|}\,, \qquad\qquad \qquad 
P_s(t,s)~ \doteq~ {\partial \over\partial s}P(t,s)\,.\eeq
The above implies
\bel{PP}
P(t,s)~=~\int_0^s (1-e^{- \alpha(t-\sigma)})\bfk(t,\sigma)\, d\sigma\,.
\eeq
We shall always assume that the curvature vanishes at the tip
of the stem, so that
\bel{tip}{\partial \over\partial s}\bfk(t,s)\bigg|_{s=t}~=~0.\eeq
%If there is no response to gravity, then 
%$${\partial\over\partial t}\bfk(t,s)~=~0,$$  and any portion of 
%the stem, once created, does not change its direction.

Our description of the growing stem takes into account:
\begi

\item[(i)] the upward bending, as a response to gravity,
\item[(ii)] an additional bending, in case of a vine clinging
to branches of other plants,
\item[(iii)] the reaction produced by obstacles,
\item[(iv)] the linear elongation.
\endi
Without loss of generality, one can assume that $P(t,0)=0\in\R^3$.
Most of our analysis will be concerned with the limit case where
 $\alpha\to +\infty$, so that  
$d\ell = d\sigma$ and (\ref{PP}) simplifies to
\bel{P}
P(t,s)~=~\int_0^s \bfk(t,\sigma)\, d\sigma\,.
\eeq
As shown in Section~\ref{s:6}, all results can be extended to the case
$0<\alpha<\infty$, with only minor changes in the proofs.

\begin{figure}[ht]
\centerline{\hbox{\includegraphics[width=12cm]{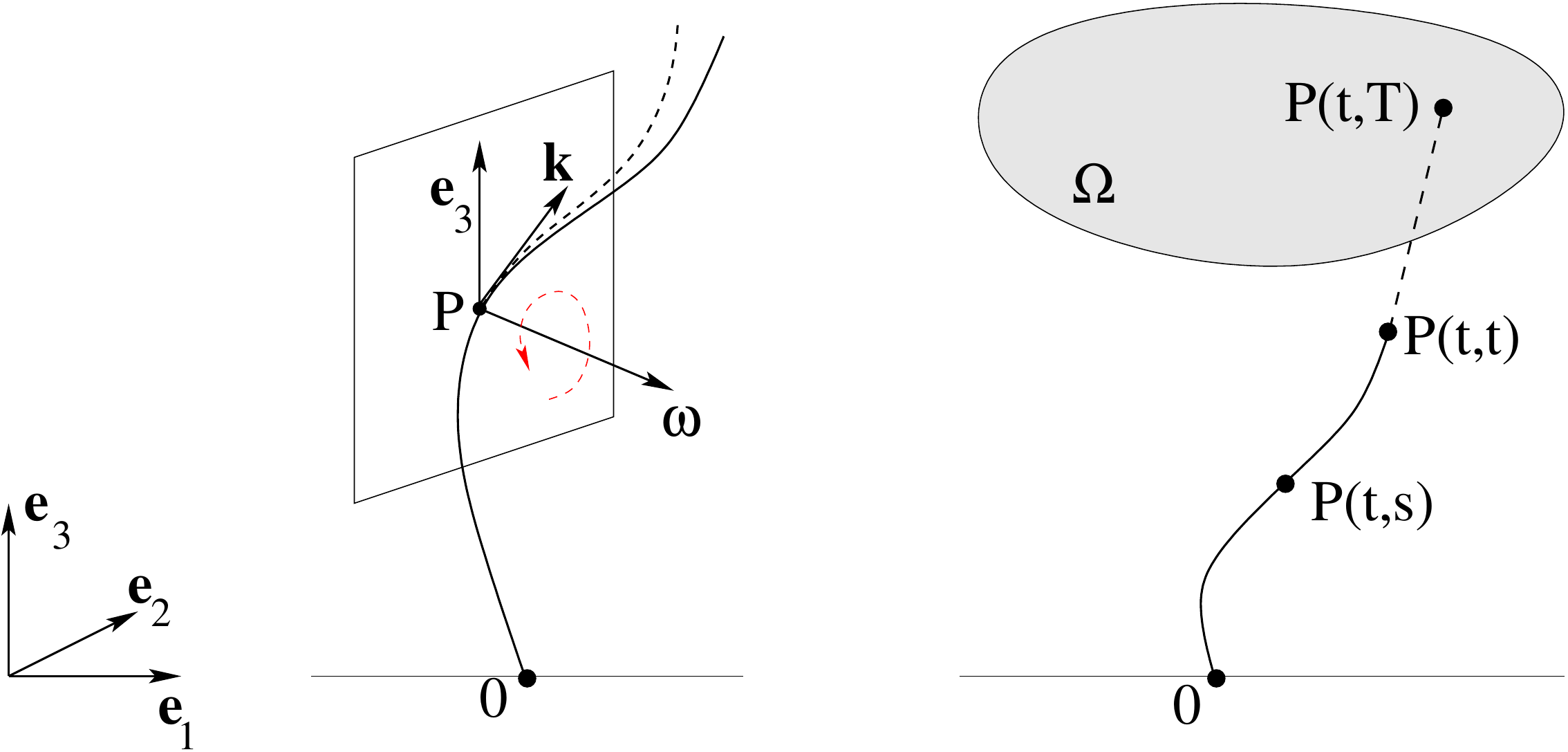}}}
\caption{\small Left: at any point $P= P(t,\sigma)$ along the stem, if the tangent 
vector $\bfk$ is not vertical, consider the plane spanned by $\bfk$ and $\bfe_3$.   Then the change in curvature of the stem at $P$ produces
a slight rotation of all points $P(t,s)$ with $s\in [\sigma,t]$.
The angular velocity is given by the vector $\omega(\sigma)$. 
Right: At a given time $t$, the curve $P(t,\cdot)$ is parameterized by $s\in [0,t]$.
It is convenient to 
prolong this curve by adding a segment of length $T-t$ at its tip
(dotted line, possibly entering inside the obstacle).
This yields an evolution equation on a 
fixed functional space $H^2([0,T];\,\R^3)$.}
\label{f:sg95}
\end{figure}

\subsection{Response to gravity.}
To model the response to gravity,
we assume that, if a portion of the  stem is not 
vertical, a local change in the curvature will be produced, affecting
the position of the upper section of the stem.

More precisely, 
let $\{\bfe_1, \bfe_2, \bfe_3\}$ be  the 
standard orthonormal basis in $\R^3$, with $\bfe_3$ oriented in the upward direction.
At every point $P(t,\sigma)$, $\sigma\in [0,t]$,
consider the cross product 
$$\omega(t,\sigma)~\doteq~\bfk(t,\sigma)\times \bfe_3.$$
The change in the position of points on  the stem,
in response to gravity, is described by (see Fig.~\ref{f:sg95})
\bel{F1}\bega{rl}
\ds{\partial\over\partial t}P(t,s)&\ds=~\int_0^s \kappa\, e^{-\beta(t-\sigma)}
 \,\bigl(\bfk(t,\sigma)\times \bfe_3\bigr)\times\bigl(P(t,s)-P(t,\sigma)\bigr)
\, d\sigma%\\[4mm]
~\doteq~ F_1(t,s)\,.\enda\eeq
Differentiating w.r.t.~$s$ one obtains
\bel{G1}\bega{rl}
\ds
{\partial\over\partial t}\bfk(t,s)&\ds=~\left(\int_0^s \kappa\, e^{-\beta(t-\sigma)}
 \,(\bfk(t,\sigma)\times \bfe_3)
\, d\sigma\right)\times \bfk(t,s)%\\[4mm]
~ \doteq~ G_1(t,s)\,.\enda\eeq
Notice that in the above integrands:
\begi
\item $\kappa> 0$ is a constant, measuring the strength of the response, while
$e^{-\beta(t-s)}$ is a {\bf stiffness factor}.   
It accounts for the fact that
older parts of the stem are more rigid and hence they bend
more slowly. 
\item $\omega(t,\sigma)=\bfk(t,\sigma)\times \bfe_3$ is an angular velocity, determined  by the response to gravity  at the point $P(t,\sigma)$.  This affects the upper portion of the stem, 
i.e.~all points $P(t,s)$ with 
$s\in [\sigma, t]$. 
%\item $(1-e^{- \alpha(t-\sigma)})d\sigma~=~d\ell$~ = arclength.
\endi
\v
\subsection{Clinging to obstacles.}
Some plants, rather than growing in the vertical direction, prefer
to curl around branches of other plants. To model this behavior,
we assume that the stem can feel the presence of an obstacle
within a distance $\delta_0$. 
This triggers a local change of the curvature, 
in the appropriate direction.

\begin{figure}[ht]
\centerline{\hbox{\includegraphics[width=16cm]{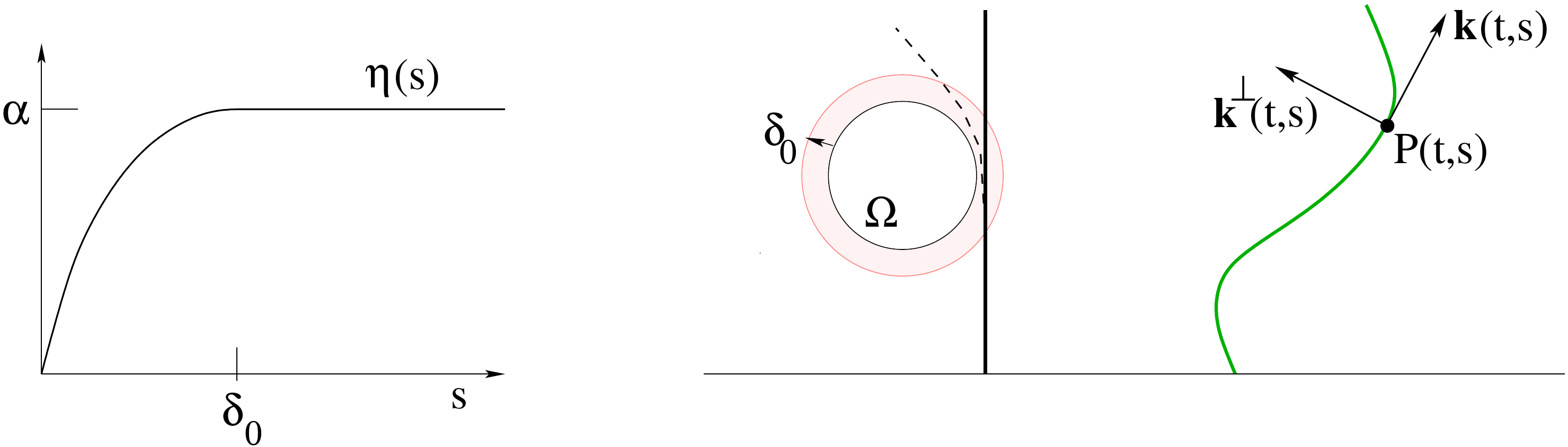}}}
\caption{\small Left: the function $\eta$ in (\ref{etaprop}).
Center: points on the stem at a distance $\leq\delta_0$  
from $\Omega$
feel the presence of the obstacle and produce an increase of curvature in the appropriate direction.}
\label{f:sg67}
\end{figure}

More precisely, let $\Omega\subset\R^3$ be a bounded open set, 
whose closure $\ov\Omega$ does not contain the origin
and  
whose boundary $\partial \Omega$ is a surface with $\C^3$ regularity.

Let 
$\eta:\R\mapsto [0,1]$ be a smooth function such that
\bel{etaprop}\eta(0) = 0,\qquad\quad
\eta(s)= \eta_0\quad \hbox{for}~~ x\geq \delta_0,\qquad\qquad \eta'\geq 0,~~ \eta''\leq 0\,,\eeq
for some constants $\eta_0, \delta_0>0$.
For $x\in \R^3\setminus\Omega$, we then set
$$\psi(x)~\doteq~\eta\bigl(d(x,\Omega)\bigr).$$
The bending of the stem around the obstacle $\Omega$ can now be described by
\bel{F2}\bega{rl}\ds
{\partial \over \partial t} 
P(t,s)&\ds =~\int_0^s e^{-\beta(t-\sigma)}\Big(
 \nabla\psi(P(t,\sigma))\times
\bfk(t,\sigma)\Big)\times\bigl(P(t,s)-P(t,\sigma)\bigr)
%\bigl(1-e^{- \alpha(t-\sigma)}\bigr)
\, d\sigma%\\[4mm]
~\doteq~F_2(t,s)\,.\enda\eeq
As before, a differentiation w.r.t.~$s$ yields
\bel{G2}\bega{rl}\ds
{\partial \over \partial t} 
\bfk(t,s)&\ds=~\left(\int_0^s e^{-\beta(t-\sigma)}\Big(
 \nabla\psi(P(t,\sigma))\times
\bfk(t,\sigma)\Big)%\bigl(1-e^{- \alpha(t-\sigma)}\bigr)
\, d\sigma\right)\times \bfk(t,s)
%\\[4mm]&
~\doteq~G_2(t,s).\enda\eeq
\v
\subsection{Unilateral constraints.}
Finally, we seek 
to model a family of  admissible reactions produced by an obstacle
$\Omega$, which guarantee that the stem will never penetrate inside 
$\Omega$.

As a preliminary, 
consider a stem which partly lies inside the open region $\Omega$.
Call 
\bel{sdist}\Phi(x)~\doteq~\left\{ \bega{rl}
 d(x,~\partial \Omega)\qquad&\hbox{if}~~x\notin\Omega,\\[3mm]
-d(x,~\partial \Omega) \qquad&\hbox{if}~~x\in\Omega,\enda\right.\eeq
 the signed distance of a point $x\in \R^3$ to the boundary $\partial \Omega$.
Let $s'\in [0,t]$ be fixed. 
If $\Phi(P(t,s'))<0$, consider the problem of  bending the stem, so that the point $P(t,s')$ is pushed out of the obstacle.   
Calling $\omega(\sigma)\in \R^3$ 
a rotation vector  at the point $P(t,\sigma)$, in first approximation
the displacement of the 
point $P(t,s')$ on the stem is computed by
\bel{disp}
\Tilde P(t,s')-P(t,s')~=~\int_0^{s'} \omega(\sigma)\times 
\bigl(P(t,s')-P(t,\sigma)\bigr)\, d\sigma\,.\eeq
We seek a function $\omega(\cdot)$ which minimizes the elastic 
deformation energy
\bel{E3}\E~\doteq~{ 1\over 2}\int_0^{s'} e^{\beta(t-\sigma)}
 |\omega(\sigma)|^2\,d\sigma\,,\eeq
subject to
\bel{paw1}\nabla\Phi(P(t,s'))\cdot \bigl(\Tilde P(t,s')-P(t,s')\bigr)+\Phi(P(t,s'))~=~0. \eeq
Here the factor $e^{\beta(t-\sigma)} $ accounts for the fact that older 
sections of the stem are stiffer, and offer more resistance to bending.
Inserting (\ref{disp}) in (\ref{paw1}) we obtain
\bel{paw2}\nabla\Phi(P(t,s'))\cdot\left( \int_0^{s'} \omega(\sigma)\times 
\bigl(P(t,s')-P(t,\sigma)\bigr)\, d\sigma\right)+\Phi(P(t,s'))~=~0. \eeq
To derive necessary conditions for optimality, consider a family of perturbations
$$\omega_\ve(\sigma)~=~\omega(\sigma)+\ve\, \tilde\omega(\sigma)\,.$$
Differentiating w.r.t.~$\ve$, at $\ve=0$ we  
obtain
\bel{nc1}  \int_0^{s'} e^{\beta(t-\sigma)} \omega(\sigma)\cdot\tilde\omega(\sigma)
\,d\sigma+\lambda\, \nabla\Phi(P(t,s'))\cdot\left( \int_0^{s'} \tilde \omega(\sigma)\times 
\bigl(P(t,s')-P(t,\sigma)\bigr)\, d\sigma\right)~=~0,\eeq
where the  constant $\lambda$ is a suitable Lagrange multiplier.

Using the property of the mixed product ${\bf a}\cdot({\bf b}\times{\bf c})=
{\bf b}\cdot({\bf c}\times {\bf a})$, one obtains
\bel{nc2}  \int_0^{s'} e^{\beta(t-\sigma)} \omega(\sigma)\cdot\tilde\omega(\sigma)
\,d\sigma+\lambda\, \left( \int_0^{s'} \tilde \omega(\sigma)
\cdot \Big(\nabla\Phi(P(t,{s'}))\times 
\bigl(P(t,{s'})-P(t,\sigma)\bigr)\Big)\, d\sigma\right)~=~0.\eeq
Since (\ref{nc2}) must hold for all perturbations $\tilde\omega(\cdot)$, 
this implies
\bel{nc3}
e^{\beta(t-\sigma)}\omega(\sigma)~=~\lambda\,\nabla\Phi(P(t,{s'}))\times 
\bigl(P(t,{s'})-P(t,\sigma)\bigr),\eeq
for some constant $\lambda\in\R$ and all $\sigma\in [0,{s'}]$.
Imposing the boundary condition (\ref{paw1}), we conclude
\bel{nc4}
\omega(\sigma)~=~\lambda\, e^{-\beta(t-\sigma)}\nabla\Phi(P(t,{s'}))\times 
\bigl(P(t,{s'})-P(t,\sigma)\bigr),\eeq
where the constant $\lambda$ is determined by the identity (\ref{paw1}).
Namely
\bel{paw4}\bega{l}
\ds\nabla\Phi(P(t,{s'}))\cdot\left( \int_0^{s'} \lambda e^{-\beta(t-\sigma)}
\Big(\nabla\Phi(P(t,{s'}))\times 
\bigl(P(t,{s'})-P(t,\sigma)\bigr)\Big)\times 
\bigl(P(t,{s'})-P(t,\sigma)\bigr)\, d\sigma\right)\\[4mm]
\qquad\qquad\ds +\Phi(P(t,{s'}))~=~0. \enda\eeq
Using the vector identities
\bel{pid1}({\bf a}\times{\bf b})\times{\bf c}~=~
({\bf c} \cdot {\bf a}){\bf b}-({\bf c} \cdot {\bf b}){\bf a}  \,,\eeq
\bel{pid2}{\bf b}\cdot \Big( (\bfb\times\bfc)\times\bfc\Big)~
=~(\bfb\cdot \bfc)^2 - |\bfb|^2|\bfc|^2,\eeq
and recalling that $\Phi(P(t, s'))<0$, from (\ref{paw4}) we obtain
\bel{kap}\bega{rl}
\lambda&\ds=~\Bigg[ \int_0^{s'} e^{-\beta(t-\sigma)} \bigg\{
\Big| \nabla\Phi(P(t,{s'}))\Big|^2\,
\Big|P(t,{s'})-P(t,\sigma)\Big|^2\\[4mm]
&\qquad\qquad \ds -
\Big( \nabla\Phi(P(t,{s'}))\cdot
\bigl(P(t,{s'})-P(t,\sigma)\bigr)\Big)^2\bigg\}\,d\sigma
\Bigg]^{-1}\Phi((P(t,{s'}))~ \leq~ 0.\enda\eeq
Notice that the integral in (\ref{kap}) vanishes only if
the vector $\nabla \Phi(P(t,{s'}))$ is parallel to all vectors
$P(t,{s'})-P(t,\sigma)$.  Examples of these ``bad" configurations 
are  shown in Figure~\ref{f:sg55}, right.
%\begin{figure}[ht]
%\centerline{\hbox{\includegraphics[width=5cm]{FIG/sg52.eps}}}
%\caption{\small An admissible reaction $v(\cdot)\in \Gamma$, 
%produced by the obstacle $\Omega$.}
%\label{f:sg52}
%\end{figure}

Next, assume that at a given time $t$ the stem lies entirely outside the obstacle, but 
part of it touches the boundary. Call
\bel{chit}
\chi(t)~\doteq~\Big\{ s\in [0,t]\,;~~
P(t,s)\in ~\partial\Omega\Big\}
\eeq
 the set where the stem touches the obstacle.
For $s\in \chi(t)$, let $\bfn(t,s)$ be the unit outer normal to the 
boundary $\partial\Omega$ at the point $P(t,s)$.  

Motivated by the previous analysis,
we define the {\bf cone of admissible velocities} produced by the 
obstacle reaction
to be the set of velocity fields
\bel{RC} \bega{l}
\Gamma(t)~\doteq~\ds \Bigg\{ \bfv:[0,t]\mapsto \R^3\,;
~~\hbox{there exists a  positive measure $\mu$ supported on $\chi(t)$
such that}\\[4mm]
\ds  \bfv(s)=-\int_0^s e^{-\beta(t-\sigma)}\left( \int_\sigma^t  
\Big(\bfn(t,s')\times 
\bigl(P(t,s')-P(t,\sigma)\bigr)\Big)   d\mu(s') \right) \times 
\bigl(P(t,s)-P(t,\sigma)\bigr)\, d\sigma %\hbox{ for all $\zeta\in [0,t]$}
\Bigg\}.
\enda
\eeq
Note that, for every $s'\in \chi(t)$, the reaction produced by the obstacle
can yield a deformation of the stem  given by
$$P_t(t,s)~=~\int_0^s \omega(\sigma)\times (P(t,s)-P(t,\sigma))\, 
d\sigma,$$
where, by (\ref{nc4}), 
$$\omega (\sigma)~=~\lambda(s')\cdot e^{-\beta(t-\sigma)}
\bfn(t,s')\times 
\bigl(P(t,s')-P(t,\sigma)\bigr),$$
for some $\lambda(s')\leq 0$.
Integrating over all points $s'\in \chi(t)$, with arbitrary choices of 
the factor $\lambda(s')\leq 0$, we obtain (\ref{RC}).
\v
{\bf Remark 1.}
One may adopt a more accurate model, where the 
bending and twisting of the stem are penalized in different ways.
More precisely, at each point $P(t,\sigma)$ one may split
the rotation vector into a component parallel to $\bfk(t,\sigma)$
(twisting) and a component perpendicular to $\bfk(t,\sigma)$
(bending):
$$\omega(\sigma)~=~\omega^{twist}(\sigma) +\omega^{bend}(\sigma),$$
where
$$\omega^{twist}(\sigma)~=~\Big( \bfk(t,\sigma)\cdot \omega(\sigma)\Big)\,\bfk
(t,\sigma),\qquad\qquad \omega^{bend}(\sigma)
 ~=~\omega(\sigma)-\omega^{twist}(\sigma).$$
The energy functional $\E$ in (\ref{E3}) can then be replaced by
\bel{E33}\E~\doteq~{ 1\over 2}\int_0^{s'} e^{\beta(t-\sigma)}
 \Big(c_1\,\bigl|\omega^{twist}(\sigma)\bigr|^2+c_2
 \bigl|\omega^{bend}(\sigma)\bigr|^2
 \Big)\,d\sigma\,,\eeq
for suitable constants $c_1,c_2$. 
When $c_1=c_2$, this is equivalent to (\ref{E3}).

\v
\subsection{Summary of the equations.}
Taking into account all terms (i)--(iii), the evolution of the stem in the presence of obstacles can be described by 
\bel{E1}
P_t(t,s)~=~F_1(t,s)+F_2(t,s) +\bfv(t,s).\eeq
Here $F_1,F_2$ are the integral terms defined 
at (\ref{F1}), (\ref{F2}), 
while
%\bel{F3}
%F_3(t,s)~\doteq~\int_0^t \alpha e^{-\alpha(t-\sigma)} \bfk(t,\sigma)\, d\sigma\eeq
%accounts for linear elongation.  Moreover,
$\bfv(t,\cdot)\in \Gamma(t)$ is an admissible reaction, 
in the cone defined at 
(\ref{RC}).
The equation (\ref{E1}) needs to be solved on a domain 
of the form 
\bel{Ddef}\D~\doteq~\bigl\{(t,s)\,;~~t\geq t_0\,,~~s\in [0,t]\bigr\}, \eeq 
with initial and boundary conditions
\bel{iP}
P(t_0,s)~=~\ov P(s),\qquad \qquad s\in [0, t_0],\eeq
\bel{bP}
P_{ss}(t,s)\bigg|_{s=t}~=~0,\qquad\qquad t> t_0\,,\eeq
and the  constraint
\bel{Pout}
P(t,s)~\notin~\Omega\qquad\forall ~(t,s)\in \D.\eeq
Recalling (\ref{P}), one obtains an equivalent evolution equation for the
unit tangent vector $\bfk$, namely 
\bel{E2}
\bfk_t(t,s)~=~G_1(t,s)+ G_2(t,s) +\bfh(t,s).\eeq
Here $G_1,G_2$ are the integral terms defined at (\ref{G1}), (\ref{G2}), respectively.
Moreover,
$\bfh(t,\cdot)$ 
is any element of the admissible cone
\bel{RC'}\bega{rl}
\Gamma'(t)&\doteq~\ds \Bigg\{ \bfh:[0,t]\mapsto \R^3\,;~~
\hbox{there exists a  positive measure $\mu$ supported on $\chi(t)$
such that}\\[4mm]
&\ds  \bfh(s)~=~ - \int_0^{s}\left(\int_\sigma^t  e^{-\beta(t-\sigma)}
\bfn(t,s')\times 
\bigl(P(t,s')-P(t,\sigma)\bigr)\,d\mu(s')\right) d\sigma\times 
\bfk(t,s)\Bigg\}.
\enda
\eeq
The equation (\ref{E2}) should be solved on the domain $\D$ in (\ref{Ddef}), 
with initial and boundary conditions
\bel{ik}\bfk(t_0, s)~=~\ov \bfk(s),\qquad \qquad s\in [0, t_0]\,,\eeq
\bel{bk}\bfk_s(t,s)\bigg|_{s=t}~=~0\qquad\qquad t>t_0\,.\eeq
\v
\subsection{The two-dimensional case}
\label{s:35}
In the planar case $n=2$, the evolution equation for the growing stem 
takes a simpler form.
For any vector $\bfv=(v_1,v_2)$, let $\bfv^\perp= (-v_2, v_1)$ be the  perpendicular 
vector obtained by a counterclockwise rotation of $\pi/2$.  Setting 
$\bfk=(k_1,k_2)$,
% and denoting by $\langle\cdot, \cdot \rangle$ the Euclidean inner product,
the equations (\ref{E2}) can be written as 
\bel{2d}\bega{rl}\ds
\bfk_t(t,s)&\ds=~\left(\int_0^s \kappa\, e^{-\beta(t-\sigma)}
  k_1(t,\sigma) \, d\sigma\right)\bfk^\perp(t,s)\\[4mm]
&\ds\qquad -\left(\int_0^s e^{-\beta(t-\sigma)}
\Big( \nabla \psi(P(t,\sigma))\cdot \bfk^\perp(t,\sigma)\Big)
\, d\sigma\right)\bfk^\perp(t,s)\\[4mm]
&\ds\qquad -
\Bigg(\int_0^s \left(\int_\sigma^t 
e^{-\beta(t-\sigma)} \Big( \bfn(t, s')\cdot \bigl(P(t, s')-P(t,\sigma)
\bigr)^\perp
\Big) d\mu(s')\right) d\sigma\Bigg)\bfk^\perp(t,s).
\enda
\eeq
Here $\mu$ is any positive measure supported on the 
contact set $\chi(t)$ in (\ref{chit}). As before, 
for $s'\in \chi(t)$ we denote by 
$\bfn(t, s')$ the unit outer normal to $\Omega$ at the 
boundary point  $P(t,s')\in \partial \Omega$.
\v

\section{Statement of the main results}
\label{s:2}
\setcounter{equation}{0}
At each time $t$, the position of the stem is described by a map
$P(t,\cdot)$ from $[0,t]$ into $\R^3$.    Of course, the domain of this map grows with
time.   
It is convenient to reformulate our model as an evolution problem 
on a functional space independent of $t$.   For this purpose, we fix $T>t_0$
and consider the Hilbert-Sobolev space $H^2([0,T];\, \R^3)$.  
Any function $P(t,\cdot)\in H^2([0,t];\, \R^3)$ will be canonically extended to 
$H^2([0,T];\, \R^3)$ by setting (see Fig.~\ref{f:sg95}, right)
\bel{Pex}
P(t,s)~\doteq~P(t,t) + (s-t)P_s(t,t)\qquad\qquad \hbox{for}~~s\in [t,T]\,.\eeq
Notice that the above extension is well defined because 
$P(t, \cdot)$ and $P_s(t,\cdot)$
are continuous functions. In all of the following analysis, we
shall study
functions defined on a domain of the form
\bel{DT}\D_T~\doteq~\bigl\{(t,s)\,;~~0\leq s\leq t,~~
t\in [t_0,T]\bigr\}, \eeq
and extended to the rectangle $[t_0,T]\times [0,T]$ as in (\ref{Pex}).
In particular, the partial derivative $P_s(t,s)$ will be constant 
for $s\in [t,T]$. This will 
already account for the boundary condition (\ref{bP}).

Adopting the notation $a\wedge b\doteq \min\{a,b\}$,
we thus consider an  evolution problem on the space $H^2([0,T];\R^3)$,
having the more general form
\bel{E5}
P_t(t,s)~=~\int_0^{s\wedge t} \Psi\Big(t,\sigma, P(t,\sigma), 
P_s(t,\sigma)\Big)\times 
\Big(P(t,s)-P(t,\sigma)\Big)
d\sigma
+\bfv(t,s).\eeq
Here $s\in [0,T]$,  $\Psi: \R\times \R\times \R^3\times\R^3\mapsto \R^3$ is a smooth function, and $\bfv(t,\cdot)$ is an admissible velocity field produced by the constraint reaction.   More precisely, given the configuration
$P(t,\cdot)$ of the stem at time $t$, 
the cone of admissible velocities is defined as
\bel{RC5} \bega{l}
\Gamma(t)~\doteq~\ds \Bigg\{ \bfw:[0,T]\mapsto \R^3\,;
~~\hbox{there exists a  positive measure $\mu$, supported on} 
\\[4mm ]\qquad\qquad\qquad \hbox{ the coincidence set 
$\chi(t)$ in (\ref{chit}),
such that for every $s\in [0, T]$ one has}\\[4mm]
\ds  \bfw(s)=-\int_0^s e^{-\beta(t-\sigma)}\left( \int_\sigma^t  
\Big(\bfn(t,s')\times 
\bigl(P(t,s')-P(t,\sigma)\bigr)\Big)   d\mu(s') \right) \times 
\bigl(P(t,s)-P(t,\sigma)\bigr)\, d\sigma 
\Bigg\}.
\enda
\eeq
%Notice that the right hand side of (\ref{E5}) is always perpendicular
%to the tangent vector $\bfk(t,s)\doteq P_s(t,s)$. 
%As a consequence, the identities
%$$|\bfk(t,s)|~ =~ |P_s(t,s)|~ =~ 1$$
%remain always valid, provided they hold at the initial time $t=t_0$. 

{\bf Remark 2.} In view of (\ref{F1}) and (\ref{F2}), we can write 
the evolution equation (\ref{E1}) in the form (\ref{E5}) by taking
\bel{PSD}
\Psi(t,\sigma, P, \bfk)~\doteq~e^{-\beta(t-\sigma)}\Big( \kappa(\bfk\times \bfe_3)
+ (\nabla\psi(P)\times \bfk)\Big).\eeq

Before stating our main
existence theorem, we introduce a precise definition of solution.
\v
{\bf Definition 1.}  {\it 
We say that a function $P=P(t,s)$, defined for 
$(t,s)\in [t_0, T]\times [0,T]$ is a solution to the equation (\ref{E5})-(\ref{RC5})
with initial and boundary  conditions (\ref{iP})--(\ref{Pout}) if the following holds.
\begi
\item[(i)] The map
$t\mapsto P(t,\cdot)$ is Lipschitz continuous from $[t_0,T]$ into $H^2([0,T];\,\R^3)$.
\item[(ii)]  For every $t,s$ one has
\bel{inte}\bega{rl}
P(t,s)&\ds=~P(t_0,s)+ \int_{t_0}^t \int_0^{s\wedge t} \Psi\Big(\tau,\sigma, P(\tau,\sigma), 
P_s(\tau,\sigma)\Big)\times 
\Big(P(\tau,s)-P(\tau,\sigma)\Big)
d\sigma\, d\tau\\[4mm]
&\qquad\quad \ds+\int_0^t\bfv(\tau,s)\, d\tau\,,\enda\eeq
where each $\bfv(\tau,\cdot)$ is an element of the cone $\Gamma(\tau)$ defined as in 
(\ref{RC5}).
\item[(iii)] The initial conditions hold:
\bel{ic4}
P(t_0,s)~=~\left\{\bega{cl}\ov P(s)&\qquad\hbox{if}~~s\in [0, t_0],\\[4mm]
\ov P(t_0)+(s-t_0)\ov P\,'(t_0)&\qquad\hbox{if}~~s\in [t_0,T].
\enda\right.\eeq
\item[(iv)]
The pointwise constraints hold:
\bel{pc4}
 P(t,s)~\notin~\Omega\qquad\quad\forall t\in [t_0,T],~~s\in [0,t].\eeq
 \endi
}

\begin{figure}[ht]
\centerline{\hbox{\includegraphics[width=14cm]{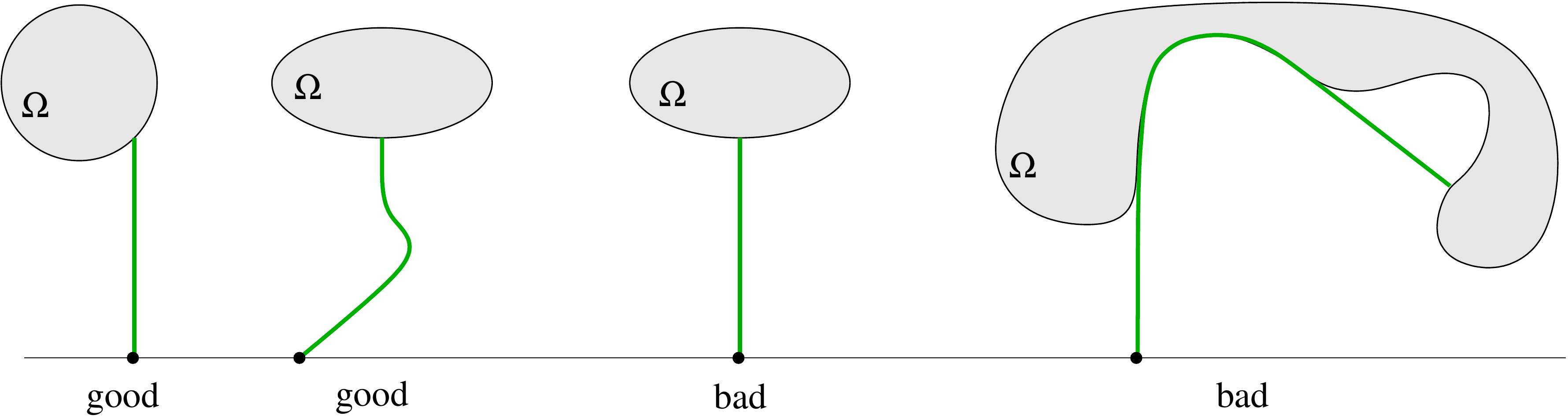}}}
\caption{\small For the two initial configurations on the 
left,  the constrained growth equation 
(\ref{E1}) admits a unique solution. %On the other hand,
The two configurations on the right satisfy the condition
{\bf (B)}. In such cases, the Cauchy problem is ill posed. }
\label{f:sg55}
\end{figure}

Given an initial data $P(t_0,s)=\ov P(s)$, 
our main result states local existence of solutions, locally in time,
as long as the following breakdown configuration
is not attained
(see Fig.~\ref{f:sg55}).
\begi
\item[{\bf (B)}] {\it The tip of the stem touches the obstacle perpendicularly, namely
\bel{bad1} \ov P(t_0)~\in~ \partial \Omega\,,\qquad\qquad 
\ov P_s(t_0)~=~-\bfn(\ov P(t_0)).\eeq 
Moreover, 
\bel{bad2} \ov P_{ss}(s)~= 0\qquad  
\hbox{for all $ s\in ~]0,t[\,$ such that~}\ov P(s)\notin \partial\Omega\,.\eeq
}
\endi
Here  $\bfn(P)$ denotes the unit outer normal to $\Omega$ at
a boundary point $P\in \partial\Omega$.  

Our  main result shows that the equations of growth 
with obstacle admit a solution (in the sense of Definition~1).
Moreover, this solution can be prolonged in time until
a breakdown configuration is  reached, as described in {\bf (B)}.
\v
{\bf Theorem 1 (existence of solutions).} 
{\it Let $\Psi$ in (\ref{E5}) be a $\C^2$ function,
and let
$\Omega\subset\R^3$ be a bounded open set with $\C^2$ boundary. 
At time $t_0$, consider the initial data (\ref{iP}), where
the curve $s\mapsto \ov P(s)$ is in 
$H^2([0,t_0];\,\R^3)$ and satisfies
\bel{ovP} \ov P(0)=0\notin \partial\Omega,\qquad\qquad 
|\ov P\,'(s)|\equiv 1,\qquad \ov P(s)\notin\Omega\qquad\forall s\in [0, t_0].\eeq
Moreover, 
assume that the condition {\bf (B)} does NOT hold.

Then there exists $T>t_0$ such that  
the equations (\ref{E5})-(\ref{RC5})
with initial and boundary  conditions (\ref{iP})--(\ref{Pout})
admit at least one solution for $t\in [t_0,T]$.

Either (i) the solution is globally defined for all
times $t\geq t_0$, or (ii) The solution can be extended 
a maximal time interval $[0,T]$, where $P(T,\cdot)$ satisfies all
conditions in {\bf (B)}.} 
\v
{\bf Remark 3.} From an abstract point of view, our evolution problem has the form 
\bel{DI}
{d\over dt} u(t)~\in~\Psi(u(t)) + \Gamma(u(t)),\eeq
where $u(t)=P(t,\cdot)\in \A\subset H^2([0,T];\,\R^3)$.  
 Here $\A$ is the set of admissible configurations, satisfying the 
 constraint (\ref{pc4}), while $\Gamma(u)$ is the cone of admissible
 velocities, defined as in (\ref{RC5}).

If $\Gamma(u)$ were the (inward pointing) normal cone to $\A$ at $u$, 
then for any two solutions $u_1,u_2$ one could expect an estimate of the type
\bel{Gro}
{d\over dt} \|u_1(t)-u_2(t)\|_{H^2}~\leq~C\cdot\|u_1(t)-u_2(t)\|_{H^2}
\,.\eeq
By Gronwall's lemma, this would imply the uniqueness and continuous dependence of solutions on the initial data \cite{CG, CMM, M, RS}.
 
\begin{figure}[ht]
\centerline{\hbox{\includegraphics[width=15cm]{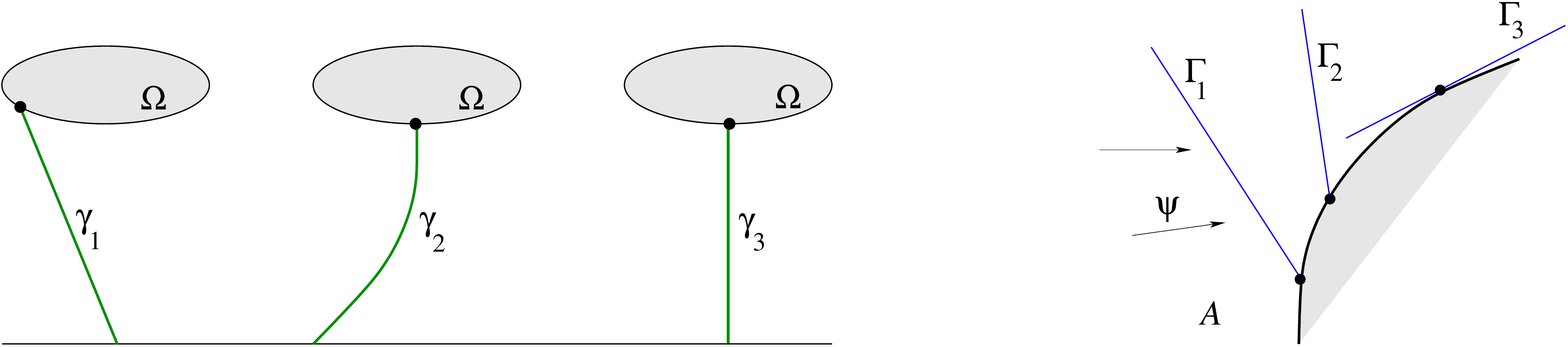}}}
\caption{\small Right: the abstract evolution equation (\ref{DI}).
In general, the cone $\Gamma$ of constraint reactions is not perpendicular to the boundary of the set $\A$ of admissible configurations.
For a stem $\gamma_3$ which satisfies all conditions in {\bf(B)}, 
the corresponding cone $\Gamma_3$ is tangent to 
$\partial \A$.}
\label{f:sg98}
\end{figure}

Unfortunately, in the present setting the 
cone $\Gamma(u)$ determined by the constraint reaction is not at all
perpendicular to the boundary of the admissible set (see 
Fig.~\ref{f:sg98}).   
For this reason, the well-posedness of the Cauchy problem 
for the growing stem with obstacles is a delicate issue, 
which will be separately addressed in the forthcoming paper \cite{BP}.

\v

\section{A  push-out operator} 
 \label{s:44}
\setcounter{equation}{0}
Consider a curve  $\gamma_0\in H^2([0,t_0];\R^3)$, parameterized by arc-length. More precisely, assume that 
\bel{gop}\gamma_0(0)~=~0\in\R^3,
\qquad |\gamma_0'(s)| ~=~ 1 ,
\qquad\gamma_0(s)\notin\Omega\quad\qquad\forall~s\in [0, t_0].\eeq
Moreover, assume that not all of the following  conditions hold:
\bel{bad5}
\gamma_0(t_0)\in \partial\Omega,\qquad\qquad \gamma'_0(t_0)= - \bfn
(\gamma_0(t_0)),\eeq
\bel{bad6} \gamma_0''(s)~=~0\qquad\hbox{for all ~$s\in [0, t_0]$~
such that }~\gamma_0(s)\notin \partial\Omega.\eeq

Given  $T>t_0$  
we can extend $\gamma_0$ to a map $[0, T]\mapsto\R^3$ by setting
\bel{extg}\gamma_0(s) ~\doteq~\gamma_0(t_0) +\gamma_0'(t_0)\, (s-t_0)
\qquad \qquad\forall s\in [t_0, T].\eeq
For a fixed radius $\rho>0$ and $T>t_0$, consider the tube 
$\V_{\rho}\subset H^2([0,T];\,\R^3)$ around
$\gamma_0$, defined by
\bel{tube}\bega{l}
\V_\rho~\doteq~\bigg\{ \gamma\in H^2([0,T];\,\R^3),~~\gamma(0)=0,\quad \gamma'(0)
= \gamma'_0(0), \\[4mm]
\qquad \qquad\qquad \ds\qquad|\gamma'(s)|=1\quad\forall s\in [0,t],
\qquad \int_0^t |\gamma''(s)- \gamma''_0(s)|^2\,ds\leq\rho
\bigg\}.\enda\eeq 
Given a curve $\gamma\in H^2([0,T];\,\R^3)$ and a  function 
$\omega\in \L^2([0,T];\,\R^3)$, we define the rotated curve
$\gamma_\omega$ as 
\bel{gomdef}\gamma_{\omega}(s)~\doteq~ \gamma(s)+
\int_0^s\omega(\sigma)\times  
\bigl(\gamma(s)-
\gamma(\sigma)\bigr)\,d\sigma.\eeq
Notice that, in general, the map $s\mapsto \gamma_\omega(s)$
is not an arc-length parameterization of $\gamma_\omega$.
However, for $\|\omega\|_{\L^2}$ small, 
we have
\bel{gom'}\bega{rl}
\bigl|\gamma_\omega'(s)\bigr|^2&\ds=~\left\langle 
\gamma'(s)+
\int_0^s\omega(\sigma)\, d\sigma\times  \gamma'(s)~,~
\gamma'(s)+
\int_0^s\omega(\sigma)\, d\sigma\times  \gamma'(s)\right\rangle
\\[4mm]
&\ds=~1+\left|
\int_0^s\omega(\sigma)\, d\sigma\times  \gamma'(s)\right|^2~=~1+
\O(1)\cdot\|\omega\|^2_{\L^2}\,.
\enda
\eeq
Next, consider a curve $\gamma\in \V_\rho$ and a subinterval
$[0,t]\subseteq [0,T]$.   If $\gamma(s)\in\Omega$ for
some $s\in [0,t]$,  we wish to push 
$\gamma$ outside $\Omega$,
but bending the curve as little as possible.  This leads to
the constrained optimization problem
\bel{mino}
\hbox{minimize:}\qquad J(\omega)
~\doteq~\int_0^T e^{\beta(t-s)}|\omega(s)|^2\, ds,\eeq
\bel{const3}
\hbox{subject to:}\qquad \gamma_\omega(s)\notin
\Omega\qquad\forall s\in [0,t].\eeq
Notice that here we allow $\gamma_\omega(s)\in
\Omega$ for $t< s\leq T$.
Before proving the existence of a minimizer, we prove  that there
exists at least one angular velocity $\omega$ that pushes every point 
of the curve  $\gamma$ away from the obstacle.  This will be achieved by
the first two lemmas.
 In the following, $\langle\cdot,\cdot\rangle$ 
 is the Euclidean inner product in $\R^3$, while 
 $\Phi(x)$ denotes the signed distance of 
 point $x$ to $\partial\Omega$, as in (\ref{sdist}).
\v
{\bf Lemma 1.} {\it Let $\gamma_0$ be as in (\ref{gop}).  
Then there exist
$T>t_0$, $\rho>0$, and a constant $C_0$ such that the following holds.

Extend $\gamma_0$ as in (\ref{extg}) and let
$\gamma\in \V_\rho$ as in (\ref{tube}). 
Assume  $\bfv
\in H^2([0,T];\,\R^3)$, with 
\bel{53}\bfv
(0)=\bfv
'(0)=0,
\qquad\quad \bigl\langle \bfv
'(s),\, \gamma'(s)\bigr\rangle ~=~0\qquad \forall 
 s\in [0,T].\eeq 
Then there exist a unique angular velocity field   $\omega\in 
\L^2([0,T];\,\R^3)$
such that
\bel{omno}
\|\omega\|_{\L^2}~\leq~C_0~\|\bfv
\|_{H^2}\,,\eeq
\bel{perp}
\bigl\langle \omega(s),\, \gamma'(s)\bigr\rangle~ =~0\qquad\quad \forall s\in [0,T],
\eeq
\bel{54}
\bfv
(s)~=~\int_0^s \omega(\sigma)\times (\gamma(s)-\gamma(\sigma))\,d\sigma\qquad\quad \forall s\in [0,T].\eeq
}
\v
{\bf Proof.} {\bf 1.} By the initial conditions in (\ref{53}),
differentiating (\ref{54})  we see that $\omega(\cdot)$ satisfies
(\ref{54}) if and only if
\bel{55} \bfv
'(s)~=~\left(\int_0^s\omega(\sigma)\,d\sigma
\right)\times \gamma'(s)\qquad\qquad\forall s\in [0,T].\eeq
\v
{\bf 2.} 
Next, consider a family of orthonormal frames
$\{\bfe_1(s), \bfe_2(s),\bfe_3(s)\}$, with $\bfe_1(s)=\gamma'(s)$
for all $s\in [0,T]$. 
We shall determine two scalar functions 
$\omega_2,\omega_3:[0,T]\mapsto\R$ such that the vector function
$$\omega(s)~=~\omega_2(s)\bfe_2(s)+\omega_3(s)\bfe_3(s)$$ 
satisfies (\ref{54}).

Using the orthogonality assumption in (\ref{53}), we obtain
 two scalar functions $z_2,z_3$ such that
\bel{h1}\bfv
'(s)~=~z_2(s)\bfe_2(s)+z_3(s)\bfe_3(s)~=~\int_0^s
\Big(\omega_2(\sigma)\bfe_2(\sigma)+
\omega_3(\sigma)\bfe_3(\sigma)\Big)d\sigma\times \bfe_1(s).\eeq
Projecting along $\bfe_2(s)$ we obtain
$$z_2(s)~=~\int_0^s
\la\bfe_2(\sigma)\times \bfe_1(s), \,\bfe_2(s)\ra \omega_2(\sigma)\,d
\sigma+
\int_0^s\la\bfe_3(\sigma)\times \bfe_1(s), \,\bfe_2(s)\ra \omega_3(\sigma)\, d\sigma.$$
By a property of the triple product, this is equivalent to
\bel{z222}z_2(s)~=~\int_0^s\Big[
\langle\bfe_2(\sigma),\,\bfe_3(s)\rangle\, \omega_2(\sigma) + \langle\bfe_3(\sigma),\,\bfe_3(s)\rangle \,\omega_3(\sigma)
\Big] d\sigma.\eeq
Similarly, projecting along $\bfe_3(s)$ we obtain
\bel{z333} z_3(s)~=~-\int_0^s\Big[
\langle\bfe_2(\sigma),\,\bfe_2(s)\rangle\, \omega_2(\sigma) 
+ \langle\bfe_3(\sigma),\,\bfe_2(s)\rangle \,\omega_3(\sigma)
\Big] d\sigma.\eeq
Observing that all quantities $\bfv
, z_i, \omega_i,\bfe_i$ in 
(\ref{h1}) are functions in $H^1$, we can differentiate one more time
and obtain the linear system of Volterra integral equations
\bel{VIE}
\left\{ 
\bega{rl} \omega_3(s)&\ds =~z_2'(s) - \int_0^s\Big[
\langle\bfe_2(\sigma),\,\bfe_3'(s)\rangle\, \omega_2(\sigma) + \langle\bfe_3(\sigma),\,\bfe_3'(s)\rangle \,\omega_3(\sigma)
\Big] d\sigma\,,\\[4mm]
 \omega_2(s)&\ds =~-z_3'(s) - \int_0^s\Big[
\langle\bfe_2(\sigma),\,\bfe_2'(s)\rangle\, \omega_2(\sigma) 
+ \langle\bfe_3(\sigma),\,\bfe_2'(s)\rangle \,\omega_3(\sigma)
\Big] d\sigma\,.\enda\right.\eeq
\v
{\bf 3.} The unique solution to the system (\ref{VIE}) can be
obtained by a standard fixed point argument.
Adopting vector notation, set $U=\begin{pmatrix}\omega_2\cr\omega_3
\end{pmatrix}$,
$Z=\begin{pmatrix}z_2\cr z_3\end{pmatrix}$.  Then (\ref{VIE}) can be written as
\bel{VI2}
U(s)~=~Z'(s) + \int_0^s B(s,\sigma) U(\sigma)\, d\sigma~\doteq~
\P[U](s),\eeq
where the matrix $B(s,\sigma)$ has norm
$$|B(s,\sigma)|~\leq~2|\bfe_2'(s)|+2|\bfe_3'(s)|~\doteq~b(s).$$
We claim that the operator $\P$ defined at (\ref{VI2})
is a strict contraction on the space $\L^1([0,T];\,\R^2)$
with equivalent norm
$$\|U\|~\doteq~\int_0^t \exp\left\{-4\int_0^s b(\sigma)\, d\sigma\right\}
|U(s)|\, ds.$$
Indeed, for any $U_1,U_2\in \L^1$ an integration by parts yields
$$\bega{rl}\|\P[U_1]-\P[U_2]\|
&\ds\leq~\int_0^t \exp\left\{-4\int_0^s b(\sigma)\, d\sigma\right\}
 b(s) \left( \int_0^s |U_1(\sigma)-U_2(\sigma)|\, d\sigma\right) ds
\\[4mm]
&=~\ds \int_0^t {1\over 4}\exp\left\{-4\int_0^s b(\sigma)\, d\sigma\right\} |U_1(\sigma)-U_2(\sigma)| ds
\\[4mm]
&\qquad\qquad +\ds {1\over 4}\exp\left\{-4\int_0^t b(\sigma)\, d\sigma\right\}  \left( \int_0^t |U_1(\sigma)-U_2(\sigma)|\, d\sigma\right) ds
\\[4mm]
&\leq~\ds {1\over 2}\int_0^t \exp\left\{-4\int_0^s b(\sigma)\, d\sigma\right\} |U_1(\sigma)-U_2(\sigma)| ds
\\[4mm]
&=\ds~{1\over 2}\|U_1-U_2\| ds\,.
\enda
$$ 
By the contraction mapping principle, the equation (\ref{VI2})
has a unique solution in $\L^1([0,T];\,\R^2)$.
In addition, we have
\bel{UN}
\|U\|_{\L^1}~\leq~C_1\|Z'\|_{\L^1}~\leq~C_2\|Z'\|_{\L^2}\,,\eeq
for some constants $C_1,C_2$ depending on $t$ and 
on the function $b(\cdot)$. In turn, this implies
\bel{UN2}\bega{rl}
\|U\|_{\L^2}^2&\ds\leq ~2\|Z'\|_{\L^2}^2 + 2 \int_0^t b^2(s)\left(
\int_0^s
U(\sigma)\, d\sigma\right)^2ds\\[4mm]
&\leq ~2\|Z'\|_{\L^2}^2 + 2\|b\|_{\L^2}^2 \|U\|_{\L^1}^2\\[4mm]
&\leq~\bigl(2+2\|b\|_{\L^2}^2 C_2^2\bigr)\,\|Z'\|_{\L^2}^2\,.
\enda
\eeq
\v
{\bf 4.} Returning to the original variables, we see that $\|Z'\|_{\L^2}
= \O(1)\cdot \|\bfv
\|_{H^2}$. On the other hand, 
$\|b\|_{\L^2}=\O(1)\cdot
\|\gamma\|_{H^2}$ is uniformly bounded 
as $\gamma$ ranges in $\V_\rho$. This completes the proof of (\ref{omno}).
\endproof
\v
{\bf Remark 4.} In the above proof, the uniqueness of the function 
$\omega(\cdot)$ was achieved by imposing the orthogonality condition
(\ref{perp}).
Without this restriction, infinitely many solutions are possible.
For example, if $\gamma'(s)=\bfe_0$ for all $s\in[0,t]$ and 
$\omega(\cdot)$ is a solution, then 
$\tilde\omega(s)=\omega(s) + \phi(s)\bfe_0$ is another solution, 
for any scalar function $\phi$.
\v
The next lemma shows that, if a stem is not close to
a bad configuration described in {\bf (B)}, then it can be pushed
away from the obstacle by a small rotation. 
\v
{\bf Lemma 2.} {\it
Assume that the initial curve 
$ \gamma_0:[0, t_0]\mapsto \R^3\setminus\Omega$  satisfies (\ref{gop}), but
it does NOT satisfy simultaneously all conditions in
(\ref{bad5})-(\ref{bad6}). Then there exist
$T>t_0$, $\rho,\delta >0$, and a constant $C_0$ such that the 
following holds.

Assume $t\in [t_0, T]$ and consider any curve 
$\gamma\in \V_\rho$, as in (\ref{tube}).
Then there exists $\omega:[0,t]\mapsto\R^3$, with 
\bel{ompic}
\|\omega\|_{\L^2}~\leq~C_0\,,\eeq
such that, for every $ s\in [0,t]$ such that
$\bigl|\Phi(\gamma(s))\bigr|\leq\delta$, one has
\bel{fuori1}
\left\langle
\int_0^s\omega(\sigma)\times
\bigl(\gamma(s)-\gamma(\sigma)\bigr)\, d\sigma\,,~\nabla \Phi(\gamma(s))
\right\rangle
~\geq~1.\eeq
%\item[(ii)]  If  $\gamma_0(t_0)\in \partial\Omega$, then
%bel{fuori2}
%\left\langle
%\int_0^t\bigl(\Psi(\sigma)+\omega(\sigma)\bigr)\times
%\bigl(\gamma(t)-\gamma(\sigma)\bigr)\, d\sigma+\gamma'(t)\,,~\bfn(t)\right\rangle
%~\geq~1.\eeq
%\endi
}
\v

\begin{figure}[ht]
\centerline{\hbox{\includegraphics[width=14cm]{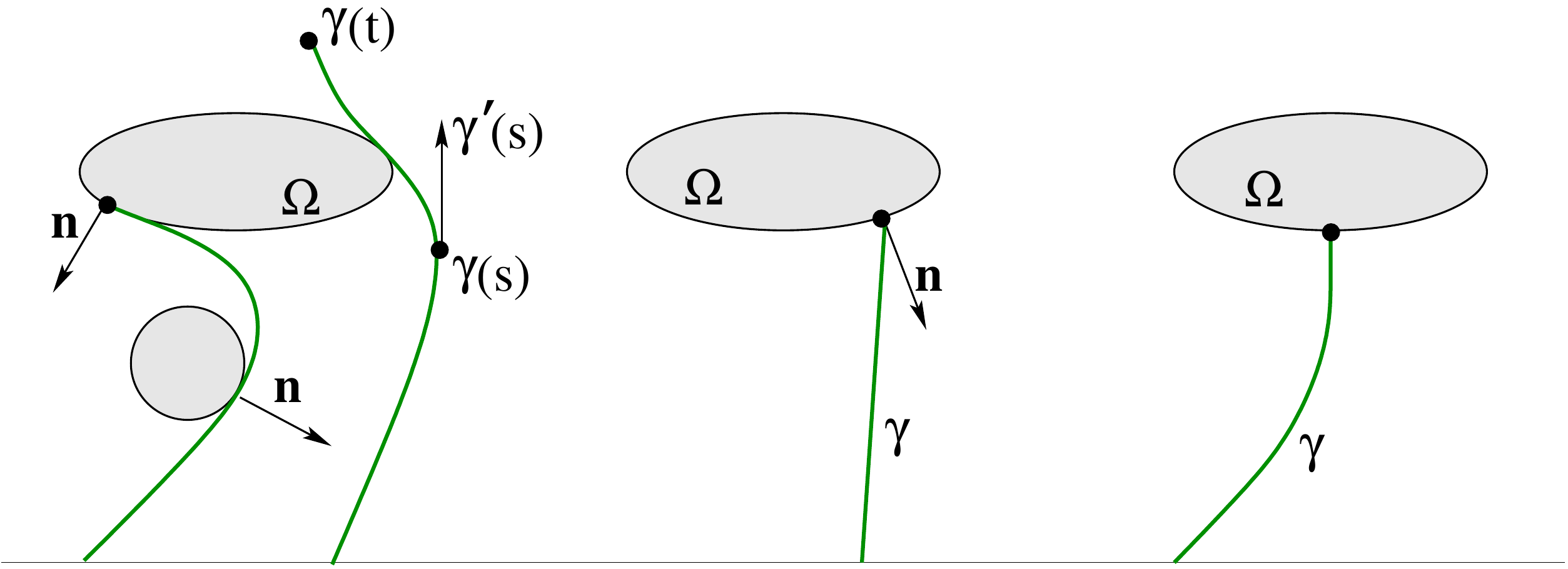}}}
\caption{\small The three cases considered in the proof of Lemma~2.}
\label{f:sg86}
\end{figure}

{\bf Proof.} Let $\bfw\in \C^2(\R^3;\,\R^3)$  be a function  which
 satisfies
\bel{vdef}
\bfw(x)~=~\left\{\bega{cl}\nabla \Phi(x)\quad&\quad \hbox{if}\quad 
|\Phi(x)|\leq
\delta_0\,,\\[4mm]
0\qquad&\quad \hbox{if}\quad 
|\Phi(x)|\geq
2\delta_0\,.\enda\right.\eeq
Since we are assuming that the boundary $\partial\Omega$ is $\C^3$,
such a function exists, provided that $\delta_0>0$ is chosen sufficiently small.
Three cases will be considered.
\v
{\bf CASE 1:}  Either $\gamma_0(t_0)\notin \partial\Omega$, or else
$\gamma(t_0)\in \partial\Omega$ and $\left\langle \gamma'(t_0), \bfn(\gamma(t_0))\right\rangle = 0$ (see Fig.~\ref{f:sg86}, left).

In this case 
%we choose $c_1>0$ large enough so that 
%\bel{fuori3}
%\left\langle c_1 \bfv_1(s)+
%\int_0^s\bigl(\Psi(\sigma)\times
%\bigl(\gamma(s)-\gamma(\sigma)\bigr)\, d\sigma\,,~\bfn(s)\right\rangle
%~\geq~0\qquad \forall ~s\in \chi_\gamma\, .\eeq
%Then 
we define
\bel{vdef1}
\tilde \bfv(s)~\doteq~ 2\bfw(\gamma(s)),\eeq
\bel{vdef2}
\bfv(s)~\doteq~\tilde \bfv(s)  -\left(\int_0^s \la \gamma'(\sigma),\,\tilde \bfv'(\sigma)
\ra\, d\sigma\right) \gamma'(s).\eeq
Observing that 
\bel{zerlo}\langle \gamma'(s),\,\bfv'(s)\rangle~=~0\qquad\qquad \forall s\in [0,t],\eeq
by Lemma~1 we can find $\omega\in \L^2([0,t])$ such that 
(\ref{omno})-(\ref{54}) hold.

To prove (\ref{fuori1}), we first observe that
\bel{v'}\left|\int_0^s\la \gamma'(\sigma),\tilde \bfv'(\sigma)\ra
d\sigma\right|
~\leq~\int_0^s |\tilde \bfv'(\sigma)|\, d\sigma
~\leq~\int_0^t \left|{d\over d\sigma} 2\bfw(\gamma(\sigma))
\right|\, d\sigma
~\leq~2t\cdot \|\bfw\|_{\C^2}\,.\eeq
We also recall that, by assumption,
$$\la \gamma_0'(s)\,,~\bfn(\gamma_0(s))\ra~=~0$$
for all $s\in [0, t_0]$ such that $\gamma_0(s)\in \partial\Omega$.
For any given $\ve_0>0$, choosing $\delta>0$ sufficiently small we 
achieve the implication
\bel{y1}\Big|\Phi(\gamma(s))\Big|~\leq~ 2\delta
\qquad\implies\qquad \Big|\la \gamma'(s)\,,~\nabla\Phi
(\gamma(s))\ra\Big|~\leq~\ve_0\eeq
for all $s\in [0,t]$.   Choosing $\ve_0= \Big(4(t_0+\delta)\|\bfw\|_{\C^2}\Big)^{-1}$, 
if $\Big|\Phi(\gamma(s))\Big|~\leq~ \delta$, by 
(\ref{v'}) and (\ref{y1}) we now have
\bel{y2}\bega{l}\ds
\left\langle
\int_0^s\omega(\sigma)\times
\bigl(\gamma(s)-\gamma(\sigma)\bigr)\, d\sigma\,,~\nabla \Phi(\gamma(s))
\right\rangle~=~\la \bfv(s)\,,~\nabla \Phi(\gamma(s))
\ra\\[4mm]
\ds
\qquad=~\la2 \nabla \Phi(\gamma(s))
\,,~\nabla \Phi(\gamma(s))
\ra - \left\langle 
\left(\int_0^s \la \gamma'(\sigma),\,\tilde \bfv'(\sigma)\ra\, d\sigma\right) \gamma'(s)\,,~\nabla \Phi(\gamma(s))
\right\rangle\\[4mm]
\ds
\qquad\geq~2-2t\|\bfw\|_{\C^2}\cdot\ve_0
~\geq~1.
\enda\eeq
\v
{\bf CASE 2:}  $\gamma_0(t_0)\in \partial\Omega$ and $-1<\langle \gamma_0'(t_0), \bfn(\gamma_0(t_0))\rangle<0$.
In other words, the tip of the stem touches the obstacle, but is neither tangent
nor perpendicular to the boundary $\partial\Omega$
(see Fig.~\ref{f:sg86}, center).

By choosing $\delta>0$ sufficiently small, we 
can find $0<\delta_1<\delta_2$ such that
\bel{d12}\Phi(\gamma(s))~\left\{\bega{l}
\leq~2\delta_0\qquad \hbox{if}~~s\in [t-\delta_1,\,t]\,,\\[3mm]
\geq~2\delta_0\qquad \hbox{if}~~s\in [t-\delta_2,\,t-\delta_1]\,,
\enda\right.\eeq
and moreover, for some $\ve_1>0$,
$$-1+2\ve_1~<~\la \gamma'(s)\,,~\nabla\Phi(\gamma(s))\ra~<~0,\qquad
\hbox{for}~~s\in [t-\delta_1,~t].$$
Consider the function 
\bel{v2}
\tilde \bfv(s)~\doteq~\left\{\bega{rl} 2\bfw(\gamma(s))
\qquad &\hbox{if} ~~s\in [0, t-\delta_1],\\[3mm]
c_2\bfw(\gamma(s))
\qquad &\hbox{if} ~~s\in [t-\delta_1,\, t],\enda\right.\eeq
where $c_2$ is a suitably large constant, whose precise value will be determined later. Since $\bfw(\gamma(s))=0$ for $s\in [t-\delta_2,
\, t-\delta_1]$, the function
$\tilde\bfv$ has the same regularity as $\gamma$. 
With this choice of $\tilde \bfv$,
we then define $\bfv(s)$ as in (\ref{vdef2}).

We claim that (\ref{fuori1}) holds.  Indeed, for $s\in [0, t-\delta_1]$
the estimates (\ref{y1})-(\ref{y2}) remain valid.
To handle the case 
 $s\in [t-\delta_1,t]$, we begin with the estimate
\bel{y5}\bega{rl}\ds
\left|\int_{t-\delta_1}^s
 \la \gamma'(\sigma),\,{d\over d\sigma}\bfw(\gamma(\sigma))
\ra\, d\sigma\right| &\ds\leq~\left| \int_{t-\delta_1}^s \la \gamma''(\sigma),\,\bfw(\gamma(\sigma))
\ra\, d\sigma\right|+\Big|\langle\gamma'(s),\, \bfw(\gamma(s))
\rangle\Big|\\[4mm]
&\ds\leq~\int_{t-\delta_1}^t \bigl|\gamma''(\sigma)\bigr|\, d\sigma + (1-2\ve_1)~\leq~1-\ve_1\,.
\enda
\eeq
We used here the fact that $\bfw(\gamma(t-\delta_1))=0$, and 
that
$\ve_1>0$ is a constant depending only on $\gamma_0$, while $\delta_1$ can be rendered as arbitrarily small by choosing
$\delta>0$ small enough.
Using the above estimate, one obtains
 \bel{y3}\bega{l}\ds
\left\langle
\int_0^s\omega(\sigma)\times
\bigl(\gamma(s)-\gamma(\sigma)\bigr)\, d\sigma\,,~\nabla \Phi(\gamma(s))
\right\rangle~=~\la \bfv(\gamma(s))\,,~\nabla \Phi(\gamma(s))
\ra\\[4mm]
\ds
\qquad=~\la c_2\,\nabla \Phi(\gamma(s))
\,,~\nabla \Phi(\gamma(s))
\ra - \left\langle 
\left(\int_0^{t-\delta_1}
 \la \gamma'(\sigma),\,\tilde \bfv'(\sigma)
\ra\, d\sigma\right) \gamma'(s)\,,~\nabla \Phi(\gamma(s))
\right\rangle\\[4mm]
\qquad\qquad \ds- \left\langle 
\left(\int_{t-\delta_1}^s
 \la \gamma'(\sigma),\,\tilde \bfv'(\sigma)
\ra\, d\sigma\right) \gamma'(s)\,,~\nabla \Phi(\gamma(s))
\right\rangle\\[4mm]
\ds
\qquad\geq~c_2-2t\|\bfw\|_{\C^2}\,\ve_0
- c_2(1-\ve_1)~\geq~1,
\enda\eeq
provided that $c_2$ was chosen sufficiently large.

As before, since (\ref{zerlo}) holds, 
by Lemma~1 we can find $\omega\in \L^2([0,t])$ such that 
(\ref{omno})-(\ref{54}) hold.
\v
{\bf CASE 3:}  $\gamma_0(t_0)\in \partial\Omega$ and $\gamma'(t_0)= - \bfn(\gamma_0(t_0))$.
Hence the tip of the stem touches the obstacle perpendicularly
(Fig.~\ref{f:sg86}, right).

Define $\bfw$  as in (\ref{vdef}).   By assumption, (\ref{bad5})
fails. 
By choosing $\delta>0$ small enough, we can thus find
 $0<a<b<t_0$ such that 
\bel{P3} \gamma'(a)~\not=~\gamma'(b),\qquad\qquad
\Phi(\gamma(s))~\geq~ 2\delta\qquad\forall s\in [a,b]\,.\eeq
In turn, we can find a smooth function $\bfz:\R\mapsto\R^3$,
supported on $[a,b]$, such that 
\bel{z1}
\int_a^b\langle \gamma'(s),\, \bfz'(s)\rangle\, ds~=~1.\eeq
Then we set
\bel{tv3}
\tilde \bfv(s)~\doteq~2\bfw(\gamma(s)) - c_3 \bfz(s),\eeq
for some constant $c_3$ large enough.
Finally, we define $\bfv(\cdot)$ as in (\ref{vdef2}).
As in the previous cases, by the identity (\ref{zerlo}) 
we can use Lemma~1 and obtain a function 
$\omega\in \L^2([0,t])$ such that 
(\ref{omno})-(\ref{54}) hold.

It remains to prove that (\ref{fuori1}) holds, provided that 
$\delta>0$ was chosen sufficiently small.   Indeed, if $\delta>0$ 
is small enough, we 
can find $0<\delta_1<\delta_2$ such that (\ref{d12}) holds
and moreover
$$-1~\leq~\la \gamma'(s)\,,~\nabla\Phi(\gamma(s))\ra~<~-{1\over 2},\qquad
\hbox{for}~~s\in [t-\delta_1,~t].$$
For $s\in [0, t-\delta_1]$, for $\delta>0$ small we
again have the bound  in (\ref{y1}). 
Hence
\bel{y8}\bega{l}\ds
\left\langle
\int_0^s\omega(\sigma)\times
\bigl(\gamma(s)-\gamma(\sigma)\bigr)\, d\sigma\,,~\nabla \Phi(\gamma(s))
\right\rangle~=~\la \bfv(\gamma(s))\,,~\nabla \Phi(\gamma(s))
\ra\\[4mm]
\ds
\qquad=~\la 2\nabla \Phi(\gamma(s))
\,,~\nabla \Phi(\gamma(s))
\ra - \left\langle 
\left(\int_0^s
 \la \gamma'(\sigma),\,\tilde \bfv'(\sigma)
\ra\, d\sigma\right) \gamma'(s)\,,~\nabla \Phi(\gamma(s))
\right\rangle\\[4mm]
\ds
\qquad\geq~2- C\,\ve_0~\geq~1,
\enda\eeq
Notice that here $C$ is a constant that depends only on $\gamma_0$, while
$\ve_0>0$ can be rendered arbitrarily small by choosing $\delta>0$
small enough.

Finally, when $s\in [t-\delta_1,\,t]$ we have
\bel{y9}\bega{l}\ds
\left\langle
\int_0^s\omega(\sigma)\times
\bigl(\gamma(s)-\gamma(\sigma)\bigr)\, d\sigma\,,~\nabla \Phi(\gamma(s))
\right\rangle~=~
\la \bfv(\gamma(s))\,,~\nabla \Phi(\gamma(s))
\ra\\[4mm]
\ds
\qquad=~\la 2\nabla \Phi(\gamma(s))
\,,~\nabla \Phi(\gamma(s))
\ra - \left\langle 
\left(\int_0^s
 \la \gamma'(\sigma),\,\tilde \bfv'(\sigma)
\ra\, d\sigma\right) \gamma'(s)\,,~\nabla \Phi(\gamma(s))
\right\rangle\\[4mm]
\ds
\qquad =~2- 
\left(\int_0^a+\int_b^s\right)
 \la \gamma'(\sigma),\,\tilde \bfv'(\sigma)
\ra\, d\sigma~\cdot \la \gamma'(s)\,,~\nabla \Phi(\gamma(s))
\ra\\[4mm]
\ds\qquad\qquad \quad -
\int_a^b
 \la \gamma'(\sigma),\,\tilde \bfv'(\sigma)
\ra\, d\sigma~\cdot \la \gamma'(s)\,,~\nabla \Phi(\gamma(s))
\ra\\[4mm]
\ds\qquad \geq~2-\left(\int_0^a+\int_b^s\right)|\tilde \bfv'(\sigma)|
\,d\sigma + c_3 \int_a^b\langle \gamma'(\sigma),\,
\bfz'(\sigma)\rangle\,d\sigma
\la \gamma'(s)\,,~\nabla \Phi(\gamma(s))
\ra\\[4mm]
\ds\qquad \geq~2-C+{c_3\over 2}\,.
\enda\eeq
We observe that here the constant $C$ depends only on $\gamma_0$,
while $c_3$ can be chosen sufficiently large so that the right hand side
of (\ref{y9}) is $\geq 1$.
This completes the proof.
\endproof
\v
{\bf Remark 5.} If all conditions (\ref{bad5})-(\ref{bad6}) hold,
then the conclusion of Lemma 2  can fail.
For example, assume that $\gamma_0$ is a segment, with
$\gamma_0(t_0)\in \partial\Omega$ and 
 $\gamma_0'(s)=\bfe_0=-\bfn(\gamma(t_0))$ for all $s\in [0, t_0]$.
Then there is no 
field of angular velocities $\omega(\cdot)$ which satisfies 
(\ref{fuori1}) at $s=t_0$.   Indeed, in this case
for every $\omega(\cdot)$ one has 
$$
\left\langle
\int_0^{t_0}\omega(\sigma)\times
\bigl(\gamma_0(t_0)-\gamma_0(\sigma)\bigr)\, d\sigma\,,~\nabla \Phi(\gamma(t_0))
\right\rangle~=~\left\langle
\int_0^{t_0}\omega(\sigma)\times
(t_0-\sigma)\bfe_0\, d\sigma\,,~-\bfe_0
\right\rangle
~=~0.$$
\v
Given a curve $\gamma:[0,T]\mapsto\R^3$ and $t\in [0,T]$, 
we introduce the quantity
\bel{Ego}
E(t,\gamma,\Omega)~\doteq~ \sup\Big\{ d(\gamma(s),\partial\Omega)\,;~~~
s\in [0,t]\,,~\gamma(s)\in\Omega\Big\},\eeq
measuring the maximum
depth at which the initial portion of $\gamma$ (i.e., 
for $s\in [0,t]$)
penetrates inside the obstacle $\Omega$.
\v
{\bf Lemma 3.} {\it For a given path 
$\gamma_0:[0, t_0]\mapsto \R^3$, assume that at least one of the conditions
in (\ref{bad5})-(\ref{bad6}) fails. 
Then there exist $T>t_0$ and $\rho>0$
such that the following holds. 

For every
 $\gamma\in \V_\rho$ and $t\in [0,T]$,  
the problem (\ref{mino})-(\ref{const3}) has a unique solution 
$\bar\omega$.
Moreover, for some constant $C_0$ one has the estimate
\bel{omes}
\|\bar \omega\|_{\L^2}~\leq~2C_0\cdot E(t,\gamma,\Omega),\eeq 
for some constant $C_0$ independent of $\gamma\in \V_\rho$.}
\v
{\bf Proof.}  {\bf 1.}
By choosing $T-t_0$ and $\rho>0$ small enough, for any
$\gamma\in \V_\rho$
an application of Lemma~2 yields the existence 
of some angular velocity $\omega(\cdot)$ such that
(\ref{ompic})-(\ref{fuori1}) hold, for some uniform constant $C_0$.
Set 
$$\hat \omega~\doteq~2E(t,\gamma,\Omega) \cdot \omega\,.$$
Then, by choosing $\delta,\rho>0$ small enough we obtain 
$$\gamma_{\hat \omega}(s)\notin\Omega\qquad\qquad\forall 
s\in [0,t]\,.$$
\v
{\bf 2.} Now
take a minimizing sequence $(\omega_n)_{n\geq 1}$.  
By the previous 
analysis, we can assume that 
\bel{onb}
\|\omega_n
\|~\leq~2C_0\cdot E(t,\gamma,\Omega) \eeq
for all $n\geq 1$.
Moreover, it is not restrictive to assume that $\omega_n(s)=0$
for $t<s\leq T$.

We then extract  a  subsequence that converges weakly
in $\L^2([0,t])$, say $\omega_n
\weakto \tilde \omega$.
Clearly
$$\|\tilde \omega\|~\leq~\liminf_{n
\to\infty} \|\omega_n
\|~
\leq~2C_0\cdot E(t,\gamma,\Omega).$$
Hence  $\tilde \omega$ achieves the minimum.
Finally, the condition
$$\gamma_{\tilde\omega}(s)\notin\Omega 
\qquad \forall ~~~s\in [0,t]$$
follows by the uniform convergence of the 
sequence $\gamma_{\omega_n
}$ on $[0,t]$.
\v
{\bf 3.}
 To prove uniqueness, consider two minimizers:
$\omega_1\not=\omega_2$, and set $\omega=(\omega_1+\omega_2)/2$.
Observe that
\bel{o12}
\|\omega\|^2~=~\|\omega_i\|^2- {1\over 4}
\|\omega_1-\omega_2\|^2\,,\qquad\qquad i=1,2\,.\eeq
If $\gamma_\omega$ lies entirely outside the obstacle, we already reach a contradiction.  In general, the definition (\ref{gomdef}) implies that
$$\gamma_\omega(s) ~=~{ \gamma_{\omega_1}(s) +\gamma_{\omega_2}(s).
\over 2}\,$$
is the midpoint of a segment with vertices outside $\Omega$.
Observe that, for every $s\in [0,t]$, one has 
$$\bigl| \gamma_{\omega_1}(s) -\gamma_{\omega_2}(s)\bigr|
~=~\O(1)\cdot \bigl(\|\omega_1\|+ \|\omega_2\|
\bigr)~ =~ \O(1) \cdot E(t,\gamma,\Omega)\,.$$
Since the boundary $\partial\Omega$ is smooth, we conclude
$$E(t,
\gamma_\omega,\Omega)~=~\O(1)\cdot E^2(\gamma,\Omega).$$
By the argument used in Lemma~3,, we can construct a perturbation
$\hat\omega$ with
$$\|\hat\omega\|~
=~\O(1)\cdot \| \gamma_{\omega_1} -\gamma_{\omega_2}\|_{\L^\infty}^2~=~
\O(1)\cdot E^2(t,\gamma,\Omega), $$
such that $\gamma_{\omega+\hat\omega}(s)\notin\Omega$ for all
$s\in [0,t]$.
We now reach a contradiction by observing that
$$\bega{l} \|\omega+\hat\omega\|^2_\beta~\leq~\|\omega\|^2\,+ 
\O(1) \cdot \|\omega\|\, \| \gamma_{\omega_1} 
-\gamma_{\omega_2}\|^2_{\L^\infty}\\[4mm]
\qquad \leq ~\|\omega_i\|^2- {1\over 4}
\|\omega_1-\omega_2\|^2 +
\O(1) \cdot \|\omega\|\, \| {\omega_1} 
-{\omega_2}\|^2~<~\|\omega_i\|^2\,,\enda$$
provided that $\|\omega\|$ is sufficiently small.
\endproof

The next lemma derives some necessary conditions for optimality,
and provides a useful representation of the minimizing function
$\bar\omega(\cdot)$.

\v
{\bf Lemma 4.}\label{Nec.Cond.} {\it 
In the same setting as Lemma~3, let $\bar\omega$ be a minimizer
for the problem (\ref{mino})-(\ref{const3}).
Then, choosing $\delta, \rho>0$
small enough, the following holds. 
For all $s\in [0,t]$ one has 
\bel{NC_omega}\bar{\omega}(s)~=~-\int_{s}^{t}\int_{[\sigma,t]}e^{-\beta(t-s)}\nabla \Phi(\gamma_{\bar{\omega}}(s'))d \mu(s')
\times \gamma'(\sigma)d\sigma ,\eeq
where $\mu$ is a positive measure, supported on the contact set 
$\chi\doteq\Big\{s\in [0,t]\,;~~\gamma_{\bar{\omega}}(s)\in \partial\Omega\Big\}
$.
}

\v
{\bf Proof.} {\bf 1.}
Calling $\Phi(\cdot)$ the signed distance from the boundary 
$\partial\Omega$, 
as in (\ref{sdist}),
the minimization (\ref{mino})-(\ref{const3}) can be reformulated as a standard problem of optimal 
control with state constraints.
Here $\omega(\cdot)$ is the control function. The state of the system is $y(s)=(y_1,y_2,y_3)(s)\in \R^{3+3+1}$, with dynamics
\bel{CS}
 \left\{\bega{rl}
 y_1'(s)&=~\omega(s),\\[3mm]
 y_2'(s)&=~ \gamma'(s)+y_1(s)\times \gamma'(s) ,\\[3mm]
 y_3'(s)&=~\frac{1}{2}e^{-\beta s}|\omega(s)|^2,\enda\right.\qquad\qquad 
 \left\{\bega{rl}
 y_1(0)&=~0,\\[3mm]
 y_2(0)&=~0,\\[3mm]
 y_3(0)&=~0,
 \enda\right.\eeq
 and the cost function is
 $$ g(y_{3}(t))=e^{\beta t}y_{3}(t).$$
Here $y_2(s)=\gamma_{\omega}(s)$. Indeed, an integration by parts yields
$$y_2(s)~=~ \gamma(s)+ \int_{0}^{s}\left(\int_0^\sigma\omega(s')\, ds' 
\right)\times \gamma'(\sigma)d\sigma ~=~\gamma(s)+\int_{0}^{s} \omega(\sigma) \times (\gamma(s)-\gamma(\sigma))d\sigma.$$
The state constraint is
\bel{scon1}
\Phi(y_2(s))~\geq~0\qquad\qquad\forall s\in [0,t].\eeq
Necessary conditions for optimality are provided 
in \cite{V}, Theorem 9.5.1. 
Namely, there exists a Lagrange multiplier  $\lambda \geq 0$,
 an absolutely continuous adjoint vector 
 $\textbf{p}(\cdot)=(p_{1},p_{2},p_{3})(\cdot)$ 
 and a non-negative Radon measure $\mu$,
 not all identically zero, such that
 the following holds.
 {\it
 \begin{itemize}
 %\item[$i)$] $ \lambda+ ||\textbf{p}||_{L^{\infty}}+||\mu||_{T.V.}=1; $
 \item[(i)]  The vector $\bfp$ provides a Carath\'eodory solution to
 the linear Cauchy problem on $[0,t]$
 \bel{ACP}
 \left\{ \bega{rl}  p'_{1}(s)&=~q_{2}(s)\times \gamma'(s),\\[3mm]
 p'_{2}(s)&=~0,\\[3mm]
 p'_{3}(s)&=0,\enda\right.
 \qquad\qquad 
 \left\{ \bega{rl}  p_{1}(t)&=~0,\\[3mm]
 q_2(t)&=~0,\\[3mm]
 p_{3}(t)&=-\lambda e^{\beta t},\enda\right.\eeq
 where 
$$
q_{2}(s)~\doteq~
\left\{ 
\begin{array}{ll}
p_{2}(s)-\int_{[0,s[}\nabla \Phi(y_{2}(\sigma))d\mu(\sigma) &\quad \mbox{if } ~0\leq s<t\,,
\\[3mm]
p_{2}(t)-\int_{[0,t]}\nabla \Phi(y_{2}(\sigma))d\mu(\sigma) &\quad  \mbox{if } ~s =t\,.
\end{array}
\right.
$$
\item[(ii)] For a.e.~$s\in [0,t]$, one has the optimality condition 
\bel{opto}
p_{1}(s)\bar \omega(s)+\frac{1}{2} e^{-\beta s} p_{3}(s)|\bar\omega(s)|^{2}~=~ \sup_{\omega \in \R^{3}} \left\{p_{1}(s) \omega+\frac{1}{2} e^{-\beta s} p_{3}(s)|\omega|^2\right\}.\eeq
 \item[(iii)] The positive 
 measure $\mu$ is supported on the region where
 the curve touches the obstacle:
 \bel{suppm}
  \mathrm{Supp}  (\mu) ~ \subseteq ~
  \bigl\{ s\in [0,t]\,;~~ \Phi(y_{2}(s))=0\bigr\}. \eeq
\end{itemize}
}
\v
{\bf 2.}
We claim that, if at least one of the conditions (\ref{mino})-(\ref{const3}) is not satisfied, 
then the above  necessary conditions hold with $\lambda=1$. 
Indeed, assume on the contrary that $\lambda=0$. Then the optimality 
condition (\ref{opto}) can be satisfied only if $p_{1}(s)\equiv 0$. In turn this implies
\bel{adj1}
p'_{1}(s)~=~\int_{[s,t]}\nabla \Phi (\gamma_{\bar{\omega}}(\sigma))d\mu(\sigma)\times \gamma'(s)~=~0\qquad \mathrm{for}\; \mathrm{all}~ s\in[0,t].
\eeq
Furthermore, since $\lambda$, $\textbf{p}$ and $\mu$ cannot all be zero, 
 this implies $\mu \neq 0$.
 
Now consider  the integral
 $$ \int_{0}^{t}\Phi(\gamma_{\bar{\omega}+\varepsilon\hat{\omega}}(\tau))d\mu(\tau),$$
where $\varepsilon>0$, $\hat\omega$ is a control which satisfies (\ref{ompic}) and (\ref{fuori1}) and $\mu$ is a measure satisfying the necessary conditions. Taking the first variation w.r.t. $\varepsilon$, one obtains
\bel{variation} \bega{l} \ds \frac{d}{d\varepsilon} \int_{0}^{t}\Phi(\gamma_{\bar{\omega}+\varepsilon\hat{\omega}}(\tau))d\mu(\tau)\Big|_{\varepsilon=0}~=~ \int_{0}^{t}\nabla \Phi(\gamma_{\bar{\omega}}(\tau))\cdot \int_{0}^{\tau}\hat\omega(\sigma)\times(\gamma(\tau)-\gamma(\sigma))d\sigma d\mu(\tau)\\[4mm]
\ds\qquad =~ - \int_{0}^{t}\hat{\omega}(\sigma)\cdot \int_{\sigma}^{t}\nabla\Phi(\gamma_{\bar{\omega}}(\tau))\times(\gamma(\tau)-\gamma(\sigma))d\mu(\tau)d\sigma 
\\[4mm]
\ds\qquad = \int_{0}^{t}\int_{0}^{\sigma}\hat\omega(\tau)d\tau \cdot \int^{t}_{\sigma} \nabla\Phi(\gamma_{\bar\omega}(\tau))d\mu(\tau)\times \gamma'(\sigma)d\sigma=\int_{0}^{t}\Bigg(\int_{0}^{\sigma}\hat\omega(\tau)d\tau \Bigg) \cdot p'_{1}(\sigma)d\sigma~=~0.
\enda\eeq
 On the other hand, integrating with respect to $\mu$ on $[0,t]$, it follows from equation (\ref{fuori1}) that 
 $$\ds \frac{d}{d\varepsilon} \int_{0}^{t}\Phi(\gamma_{\bar{\omega}+\varepsilon\hat{\omega}}(\tau))d\mu(\tau)\Big|_{\varepsilon=0}~\geq~ \mu([0,t])~>~0,$$
 which is a contradiction. This proves that $\lambda=1$.
\v
{\bf 3.} By the previous step,
the necessary conditions are satisfied with $\lambda=1$.
In particular, from the conditions (ii)-(iii) it follows that
$$ p_{1}(s)~=~-\int_{s}^{t}\int_{[\sigma,t]}\nabla \Phi(\gamma_{\bar{\omega}}(s')) d\mu(s')\times \gamma'(\sigma)d\sigma\,,$$
while
$$\bar{\omega}(s)~=~e^{-\beta(t-s)}p_{1}(s)~=~-\int_{s}^{t}\int_{[\sigma,t]}e^{-\beta(t-s)}\nabla \Phi(\gamma_{\bar{\omega}}(s')) d\mu(s')\times \gamma'(\sigma)d\sigma\,.$$
This concludes the proof.
\endproof
\v

The above construction allows us to define a nonlinear 
``push-out" operator
$\P$ as follows.   

For convenience,  given an angular velocity $\omega\in \R^3$, we shall
denote by $R[\omega]$  the $3\times 3$ rotation matrix 
\bel{R} R[\omega]~\doteq~e^A~\doteq~\sum_{k=0}^\infty {A^k\over k!},
\qquad\qquad A~
\doteq~\left(\bega{ccc}
0&-\omega_3&-\omega_2\cr
\omega_3&0&-\omega_1\cr
-\omega_2&\omega_1&0\enda\right).\eeq
Notice that, for every $\bar \bfv\in \R^3$, the image
$R^\omega\bar \bfv$ is the value at time $t=1$ of the solution to
$$\dot \bfv(t)~=~\omega\times \bfv(t),\qquad\qquad \bfv(0)=\bar \bfv.$$

Next, 
given $\gamma\in H^2([0,t];\,\R^3)$, let
$\bar \omega$ be the minimizer for the corresponding problem 
(\ref{mino})-(\ref{const3}). Using the the notation (\ref{R}), then define
\bel{pig}
\P[\gamma](s)~\doteq~\int_0^s R\left[\int_0^\sigma \bar\omega(\zeta)d\zeta\right]
\gamma'(\sigma)\, d\sigma\,.\eeq
Here $R^\omega$ is the rotation operator, defined as in (\ref{R}).
The difference is still  $\O(1)\cdot \|\bar\omega\|_{\L^2}^2$,
but it may be easier to estimate the 
difference $\|\P[\gamma]-\gamma\|_{H^2}
$.
We refer to $\P$ as the ``push-out operator".
\v
{\bf Remark 6.} In general,   $\P[\gamma]\not= \gamma_{\bar\omega}$.
Therefore,  it may not be true that 
$\P[\gamma](s)\notin\Omega$ for all $s\in [0,t]$.   However, the two curves
are very close.
Indeed,
differentiating (\ref{gomdef}) we obtain
\bel{gomp}
\gamma_{\bar \omega}'(s) ~=~\gamma'(s)+\left(\int_0^s
\omega(\sigma)\, d\sigma\right)\times \omega'(s)
~=~R\left[\int_0^s
\omega(\sigma)\, d\sigma\right] \gamma'(s)+
 \O(1)\cdot \|\bar \omega\|_{\L^1}^2\,.\eeq
 Therefore
\bel{gpg}
\Big| \gamma_{\bar\omega}(s) -\P[\gamma](s)\Big|~
\leq~\O(1)\cdot \|\bar\omega\|_{\L^1([0,t])}^2 ~=~
\O(1)\cdot \|\bar\omega\|_{\L^2([0,t])}^2\,.\eeq
In particular, by (\ref{omes}) the curve $\P[\gamma]$ can penetrate
inside the obstacle at most by the amount 
\bel{pen1}
E\Big( \P[\gamma],\,\Omega\Big)~=~\O(1)\cdot\|\bar\omega\|_{\L^2}^2~
\leq~C\cdot E^2(\gamma,\,\Omega).
\eeq
for some uniform constant $C$.
\v
\section{Proof of Theorem 1}
\label{s:4}
\setcounter{equation}{0}

{\bf 1.} Differentiating (\ref{E5}) w.r.t.~$s$, 
one obtains an  equation for the tangent vector 
$\bfk(t,\cdot)\in H^1([0,T];\,\R^3)$.   After an integration w.r.t.~$t$,
this can be written as
\bel{E6}
\bfk(t,s)~=~\bfk(t_0,s)+\int_0^t\left(\int_0^{\tau\wedge s} 
\Psi\Big(\tau,\sigma, P(\tau,\sigma), 
\bfk(\tau,\sigma)\Big)d\sigma
\right)\times\bfk(\tau,s)\, d\tau + \int_{t_0}^t
\bfh(\tau,s)\, d\tau.\eeq
Here $P(t,s)$ is given by (\ref{P}), while $\bfh(\tau,\cdot)$ 
is an element of the cone
\bel{RC7}\bega{l}
\Gamma'(\tau)~\doteq~\ds \Bigg\{ \bfh:[0,T]\mapsto \R^3\,;~~
\hbox{there exists a  positive measure $\mu$, supported on}\\[4mm]
\qquad \hbox{the contact set
$\chi(\tau)\subseteq [0,\tau]$ in (\ref{chit}), such that 
for every $s\in [0,T]$ one has}\\[4mm]
\ds \quad
 \bfh(s)~=~ - \int_0^{s}\left(\int_\sigma^\tau  e^{-\beta(\tau-\sigma)}
\bfn(\tau,s')\times 
\bigl(P(\tau,s')-P(\tau,\sigma)\bigr)\,d\mu(s')\right) d\sigma\times 
\bfk(\tau,s)\Bigg\}.
\enda
\eeq
This integral equation should be solved for $(t,s)\in [t_0,T]\times
[0,T]$, for some $T>t_0$, with initial  data
\bel{i4}
\bfk(t_0,s)~=~\bar \bfk(s)~=~\left\{\bega{rl}\ov P\,'(s)\qquad
&\hbox{if}~~s\in [0,t_0],\\[3mm]
\ov P\,'(t_0)\qquad
&\hbox{if}~~s\in [t_0,T].\enda\right.\eeq
Notice that
$\bfk_t(t,s)$ is always perpendicular to $\bfk(t,s)$. Therefore
$|\bfk(t,s)|\equiv 1$ for all $t,s$. Moreover, 
since $\chi(t)\subset [0,t]$, 
every element $\bfh(t,\cdot)
\in \Gamma'(t)$
is constant for $s\in [t,T]$.   By (\ref{E6}), this confirms that
$\bfk(t,\cdot)$ is constant for $s\in [t,T]$.  
\v
{\bf 2.} To construct a sequence of approximate solutions, we use an operator-splitting approximation scheme.
Let $\gamma_0(\cdot)=\ov P(\cdot)$ be the curve of initial data.
As in (\ref{tube}), consider a small 
neighborhood $\V_\rho$ of $\gamma_0$
for which the conclusions of Lemmas 1--3 hold.

Fix a time step $\ve>0$ and set $t_k = t_0+\ve k$.  
Assume that the approximate solution $\bfk=\bfk(t,s)$ has been constructed for all 
times $t\in [0, t_{k-1}]$ and $s\in [0,t]$. 
To extend the solution on the next time interval $[t_{k-1}, t_k]$ 
we proceed as follows.

First, we construct an approximate solution 
of the problem without obstacle. Recalling the notation
introduced at
(\ref{R}) for a rotation matrix, we define 
\bel{41}
\bfk(t,s)~=~R 
\left[(t-t_{k-1})\int_0^{s\wedge t_{k-1}} \Psi\bigl(t_{k-1},\sigma, P(t_{k-1},\sigma), \bfk(t_{k-1},\sigma)\bigr)d\sigma
\right] \bfk(t_{k-1},s), \eeq
for $t\in [t_{k-1}, t_k[\,$ and $s\in [0, T]$. 
Taking $t= t_k$, the above construction
produces a curve 
$$s~\mapsto ~\gamma(s)~=~ P(t_k-, s)~=~\int_0^s \bfk(t,\sigma)\, 
d\sigma.$$
In general, the physically meaningful portion of this curve: 
\bel{gam6}
\gamma(s)~\doteq~P(t_k-, s)\,,\qquad\qquad s\in [0, t_k],\eeq
may well 
partly lie inside the obstacle.
Using the push-out operator (\ref{pig}), we then replace 
$P(t_k-,\cdot)$ by a new curve, by setting
\bel{pig2}
P(t_{k}, s)~=~\P [P(t_{k}-,\cdot)](s)~=~\int_0^s R\left[\int_0^\sigma \bar\omega_k(\zeta)d\zeta\right]
P_s(t_k-, \sigma)\, d\sigma,\qquad s\in [0,T]
\,.\eeq
Equivalently:
\bel{k+}
\bfk(t_k,s)~=~ R\left[\int_0^s \bar\omega_k(\zeta)d\zeta\right]
\bfk(t_k-, s)\,.\eeq
As in Lemma~3, 
here $\bar \omega_k$ is the unique minimizer for the problem 
(\ref{mino})-(\ref{const3}), with $\gamma:[0, t_k]\mapsto \R^3$
as in (\ref{gam6}).
According to Lemma~4, we have the representation
\bel{omek}
\bar\omega_{k}(s)~=~-e^{-\beta(t_k-s)}\int^{t_k}_{s}
\left(\int_{[\sigma,t_k]}\nabla\Phi(\gamma_{\bar\omega_{k}}(s'))d\mu_{k}(s') \right)
\times
P_s(t_{k}-,\sigma)d\sigma,
\eeq 
where $\mu_k$ is a positive measure supported on the contact set
\bel{chik}\chi_{k}~=~\{s\in[0,t_k]:\; 
\gamma_{\bar \omega_{k}}(t,s)\in \Omega\}.\eeq
Applying Fubini's theorem, one can exchange the order of integration in (\ref{omek})
and obtain
\bel{410}\bega{rl}\ds
\omega_{k}(s)&\ds =~-e^{-\beta(t_k-s)}\int^{t_k}_{s}
\nabla\Phi(\gamma_{\bar\omega_{k}}(s'))\times \Big(\int_{s}^{s'} 
P_s(t_{k}-,\sigma)d\sigma\Big)d\mu_{k}(s')
\\[4mm]
&\ds=~-e^{-\beta(t_k-s)}\int^{t_k}_{s}\nabla\Phi(\gamma_{\bar\omega_{k}}(s'))
\times\Big(P(t_k-, s')-P(t_k-, s)\Big)\,d\mu_{k}(s').
\enda\eeq

By (\ref{omes}) and (\ref{pen1}), as long as the 
approximation remains inside the neighborhood $\V_\rho$,
there exists a constant $C_3$ such that, recalling the definition
(\ref{Ego}) one has
\bel{ol2}
\|\bar\omega_k\|_{\L^2([0, t_k])}~\leq~C_3\, 
E(t_k,\,P(t_k-,\cdot),\,\Omega)\,,\eeq
\bel{ol3} E(t_k,\,P(t_k,\cdot),\,\Omega)~\leq~C_3 
\,E^2(t_k,\,P(t_k-,\cdot),\,\Omega).\eeq
Moreover,  the depth at which the 
curve $P(t_k-,\cdot)$ penetrates the obstacle is estimated by
\bel{ol4} E\bigl(t_k,\,P(t_{k}-,\cdot),\,\Omega\bigr)~
\leq ~E(t_{k-1},\,P(t_{k-1},\cdot),
\,\Omega)+C_4\ve\,,\eeq
for some  constant $C_4$. For every $\ve>0$ sufficiently small, the above
estimates
(\ref{ol3})-(\ref{ol4})  yield the implication
$$E(t_{k-1},\,P(t_{k-1},\cdot),\,\Omega)~\leq~\ve\qquad\implies
\qquad E(t_k,\,P(t_{k},\cdot),\,\Omega)~\leq~\ve\,.$$
 By assumption, the initial condition
lies outside the obstacle, i.e.
\bel{51} E(t_k,\,P(t_0,\cdot),\,\Omega)~=~0.\eeq
By induction, for all $\ve>0$ small enough and all 
$k\geq 1$ we thus conclude
\bel{52}E(t_k,\,P(t_k,\cdot),\,\Omega)~\leq~\ve\,.\eeq
In turn, by (\ref{ol2})-(\ref{ol3}) this 
implies an estimate of the form
\bel{57} \|\bar\omega_k\|_{\L^2([0, t_k])}~\leq~C_5\,\ve\,.\eeq 
\v
{\bf 3.} 
We claim that, for every $\ve>0$ and $k\geq 1$, 
the total mass of the positive  measure $\mu_k$ in (\ref{omek})
is bounded by 
\bel{mubo}
\|\mu_{k}\|~=~\mu_k([0,t_k])~\le~C_6\,\ve,\eeq
for some  constant $C_6$, independent of $\ve,k$.
Indeed, (\ref{omes}) implies
$$|P(t_k-,s)-\gamma_{\bar\omega_k}(s)|~ =~ \O(1)\cdot
\|\bar\omega_k\|_{\L^1}~=~\O(1)\cdot\ve.$$
Taking $\gamma:[0, t_k]\mapsto\R^3$ as in (\ref{gam6}) and
integrating the left hand side of  (\ref{fuori1}) 
w.r.t.~the measure $\mu_{k}$
 we obtain
\bel{muk} \bega{l}\ds \mu_{k}([0, t_k])~\leq~
\int_{[0, t_k]}\left\langle
\int_0^s\omega(\sigma)\times
\bigl(P(t_k-, s)-P(t_k-, \sigma)\bigr)\, 
d\sigma\,,~\nabla \Phi(P(t_k-, s))
\right\rangle d\mu_{k}(s)
\\[4mm]\ds
= ~
-\int_0^{t_k}\omega(\sigma)\cdot \Big(\int_{[\sigma, t_k]}
\bigl(\nabla \Phi(P(t_k-, s))-\nabla \Phi(\gamma_{\bar{\omega}_{k}}(s))\bigr)\times \bigl(P(t_k-, s)-P(t_k-, \sigma)\bigr)
d\mu_{k}(s)\Big)d\sigma 
\\[4mm]\ds
\qquad \qquad -\int_{0}^{t_k}\omega(\sigma)\cdot 
\Big(\int_{[\sigma,t_k]}\nabla \Phi(\gamma_{\bar{\omega}_{k}}(s)) 
\times \bigl(P(t_k-, s)-P(t_k-, \sigma)\bigr)d\mu_{k}(s)\Big)d\sigma 
\\[4mm]\ds
\ds\leq~\int_0^{t_k} |\omega(\sigma)| \,C\ve \,\|\mu_k\|\, d\sigma+
\left| \int_{0}^{t}\omega(\sigma)\cdot\bar{\omega}_{k}(\sigma)d\sigma
\right|
\\[4mm]
\ds\leq~C\ve \,\|\omega(\sigma)\|_{\L^1}\,\|\mu_k\|+
\|\omega\|_{\L^2}\|\bar\omega_k\|_{\L^2}\,,
\enda\eeq
for some constant $C$ independent of $\ve,k$.
By (\ref{ompic}) we have  
$$\|\omega\|_{\L^1([0, t_k])}~\leq~\sqrt{t_k}\, \|\omega\|_{\L^2([0, t_k]}
~\leq~\sqrt T\, C_0\,,$$ while (\ref{omes}) implies
 $\|\bar\omega_k\|_{\L^1}$, $\|\bar\omega_k\|_{\L^2}=\O(1)\cdot\ve$.
Choosing $\ve>0$ small enough so that $C\ve \|\omega\|_{\L^1}<1/2$, 
from (\ref{muk}) we deduce
$${1\over 2} \,\|\mu_k\|~\leq~\|\omega\|_{\L^2}\|\bar\omega_k\|_{\L^2}
~=~\O(1)\cdot\ve.$$
This yields the estimate (\ref{mubo}).

Using (\ref{mubo}) in (\ref{410}) we can refine the estimate (\ref{57})
and conclude that
\bel{59} \|\bar\omega_k\|_{\L^\infty([0,t_k])}~\leq~C_9\ve,\eeq
for some constant $C_9$ independent of $k,\ve$. 
\v
{\bf 4.}
By the previous arguments, for every time step 
$\ve>0$ sufficiently small, we obtain
a piecewise continuous approximate solution $\bfk_\ve=\bfk_\ve(t,s)$
defined for $s\in [0,t]$ and $t\in [t_0, T_\ve]$.   Here
$T_\ve$ is the supremum of all times $\tau$ for which the corresponding curve 
$P_\ve(\tau,s)=\int_0^s\bfk_\ve(\tau,\sigma) d\sigma$ remains in 
the neighborhood $\V_\rho$.

%In this step we study the regularity of the approximations
%$\bfk_\ve$, and show that they are well defined for $t\in [t_0,T]$,
%for some $T>t_0$ independent of $\ve$.
%Toward this goal, we claim

We now observe that, for every $\ve>0$ small enough and $k\geq 1$, as long as the approximation $P_\ve(t,\cdot)$ remains inside $\V_\rho$ one has
\bel{61}
\|\bfk_\ve(t,\cdot)-\bfk_\ve(t',\cdot)\|_{H^1([0, t_{k-1}])}
~\leq~C_7\,|t-t'|\qquad\quad\forall t,t'\in [t_{k-1}, t_k[\,,\eeq
\bel{62}
\|\bfk_\ve(t_k,\cdot)-\bfk_\ve(t_k-,\cdot)\|_{H^1([0, t_k])}
~\leq~C_7\,\ve\,,\eeq
for some constant $C_7$ independent of $k,\ve$.
Indeed, the first estimate is an immediate consequence
of the boundedness of $\Psi$.
The second estimate follows from (\ref{57}).
As long as the approximations remain inside $\V_\rho$,
the bounds (\ref{61})-(\ref{62}) imply 
an estimate of the form
\bel{63}\|\bfk_\ve(t,\cdot) - \bfk_\ve (\tau,\cdot)\|_{H^1([0,\,t])}~\leq~ C_8(\ve +\tau-t)\,,\eeq 
for some constant $C_8$ and all $t<\tau$. Since by construction
$$P_\ve(t,0)=0,\qquad\qquad \bfk_\ve(t,0)~=~
\bar \bfk(0)\qquad\forall t\geq 0,$$
from (\ref{63}) we deduce 
\bel{65}P_\ve (t_k,\cdot)~\in ~\V_\rho\qquad\qquad\forall t\in [0,T],\eeq
for some $T>0$ independent of $\ve.$

Thanks to the above estimates, we conclude that
for  $\ve>0$ sufficiently small all
the approximations $\bfk_\ve=\bfk_\ve(t,s)$ are well defined.
To achieve the convergence of a subsequence, 
two observations are in order:
\begi
\item As maps from $[t_0,T]$ into $H^1([0,T])$, all functions
$t\mapsto \bfk_\ve(t,\cdot)$ have uniformly bounded total variation.
\item For every $t,\ve$, consider the difference 
of the partial derivatives
$$\bfw(t,s)~\doteq~{\partial\over\partial s}\bfk_\ve(t,s) - {\partial\over\partial s}\bar\bfk(s).$$
Then, for every $t\in [t_0,T]$, the map $s\mapsto \bfw(t,s)$
is uniformly Lipschitz continuous.  Indeed, this is because (i)
the integral in (\ref{41}) is uniformly bounded, for all $k,\ve$, and
(ii)
the function $\bar\omega_k$ in (\ref{k+}) satisfies the uniform bound 
(\ref{59}).
As a consequence, 
all functions $\bfk_\ve(t,\cdot)$ remain within a fixed
compact subset ${\mathcal K}\subset H^1([0,T];\,\R^3)$.
\endi
Thanks to the above properties, 
we can thus extract a subsequence $\ve_n\to 0$
and achieve the $H^1$ convergence 
$\bfk_{\ve_n}(t,\cdot)\to \bfk(t,\cdot)$, uniformly for
$t\in [t_0,T]$.
\v
{\bf 5.} In this step we study the convergence 
of the measures $\mu_k^\ve$.
For any fixed $\ve>0$, let $\mu^\ve$ be the 
positive measure on $[t_0,T]\times\R$ whose restriction to 
$\,]t_{k-1}, t_k]\times\R$ coincides with the 
product $\cL\otimes\mu_k^\ve$.  Here $\cL$ denotes the Lebesgue measure.
In other words, for every subinterval $[a,b]\subseteq [t_{k-1}, t_k]$
and every open set $V\subset\R$, we have
\bel{medef}\mu^\ve\bigl([a,b]\times V\bigr)~=~(b-a)\cdot\ve^{-1}\mu_k^\ve(V).\eeq
In view of the uniform bound (\ref{muk}), we can extract  a weakly convergent subsequence, so that 
$\mu^{\ve}\rightharpoonup\mu$.   More precisely (see \cite{AFP}) , there exists
a measurable family of uniformly bounded positive measures 
$\{\mu^t\,;~~t\in [0,T]\}$ such that, for every continuous function
$\vp:\R^2\mapsto\R$ one has
\bel{421}\lim_{\ve\to 0}
\int_0^T\int \vp(t,s)\,d\mu^\ve(t,s)~=~\int_0^T\int \vp(t,s)\,d\mu(t,s)=~\int_0^T 
\left(\int \vp(t,s)d\mu^t(s)\right)dt.
\eeq
For each $\ve>0$ and $t\in [0,T]$ we define
\bel{chiep}\chi_\ve(t)~\doteq~\{s\in[0,t] \,;\, \gamma_{\bar\omega_{k}}(s)\in \partial\Omega \}~\supseteq~
\mathrm{Supp}\{\mu_{k}^\ve \},\eeq
where $k$ is the unique index such that  $t\in \,]t_{k-1}, t_k]\doteq
\,]\ve (k-1),\,\ve k]$.
By (\ref{gpg}), (\ref{57}), and (\ref{63}) it follows that $\gamma_{\omega_{k}}(t,s)\rightarrow P(t,s)$ uniformly on $\mathcal{D}_{T}$. Now fix any $\delta>0$  and consider the set 
$$V^{\delta}(t)~\doteq~
\bigl\{ s\in[0,t]\,;~~ d(P(t,s), \partial\Omega)\leq \delta\bigr\}.$$ 
In view of the uniform convergence $\bfk_\ve\to \bfk$ and  
$P_\ve\to P$
as $\ve\to 0$, it follows that $\chi_\ve(t)\subseteq V^{\delta}(t)$ for all $\ve>0$ sufficiently small. 
By the weak convergence, one has 
$$\mathrm{Supp}\{ \mu^t\}~\subseteq~\limsup_{\ve\to 0}\chi_\ve(t)~\subseteq ~V^{\delta}(t)$$
for any $\delta>0$ and $t\in [0,T]$.
Since $\delta>0$ was arbitrary, this implies
\bel{423}\mathrm{supp}\bigl\{ \mu^t\}~\subseteq~\chi(t)~=~\{ s\in[0,t]\,
;~~ P(t,s)\in \partial \Omega\bigr\}\eeq
for all $t\in [0,T]$.

\v
{\bf 6.} We complete the proof by showing  that $\bfk$ provides the desired solution.
By (\ref{52}), letting $\ve\to 0$ it is clear that $P(t,s)\notin\Omega$,
for all $t\in [0,T]$ and $s\in [0,t]$.

It remains to prove that, for every $t\in [t_0,T]$,
the identity (\ref{E6}) holds, with
\bel{h4} \bega{l}\bfh(\tau,s)\ds
~=~ - \int_0^{s}\left(\int_\sigma^\tau  e^{-\beta(\tau-\sigma)}
\bfn(\tau,s')\times 
\bigl(P(\tau,s')-P(\tau,\sigma)\bigr)\,d\mu^\tau(s')\right) d\sigma\times 
\bfk(\tau,s)\\[4mm]\ds
~=~ - \int_0^{s}\left(\int_\sigma^\tau  e^{-\beta(\tau-\sigma)}
\nabla\Phi(P(\tau,s'))\times 
\bigl(P(\tau,s')-P(\tau,\sigma)\bigr)\,d\mu^\tau(s')\right) d\sigma\times 
\bfk(\tau,s).\enda\eeq
Here $\{\mu^\tau;~\tau\in [t_0,T]\}$ is the family of measures
constructed in step {\bf 5.} 

The identity (\ref{E6}) will be obtained by taking the limit as $\ve\to 0$
of the identities (\ref{41}) and (\ref{k+}) satisfied 
by $\bfk^\ve$.
By the previous analysis, as $\ve\to 0$ we have 
\bel{conv}
\bfk^\ve(t,s)~\to~\bfk(t,s),\qquad\qquad P^\ve(t,s)~\doteq~\int_0^s\bfk^\ve(t,\sigma)\, d\sigma~\to~\int_0^s\bfk(t,\sigma)\, d\sigma~=~P(t,s),\eeq
uniformly for $(t,s)\in [t_0,T]\times [0,T]$.

To handle the right hand side of (\ref{E6}),
we start with matrix estimate
$$\Big| (I+A_1+\cdots+A_n)\bfv - e^{A_1}\circ\cdots\circ e^{A_n}\bfv
\Big|~=~\O(1)\cdot\left(\sum_i|A_i|\right)^2|\bfv|,$$
which implies, as a special case:
\bel{RR}
\Big|\bfv+
(\omega_1+\cdots
+\omega_n)\times\bfv- R[\omega_1]\circ\cdots\circ R[\omega_n]\bfv\Big|
~=~\O(1)\cdot\left(\sum_i|\omega_i|\right)^2|\bfv|.
\eeq
Using the notation
$$t_k^\ve\doteq t_0 + k\ve\,,\qquad\qquad 
k_\ve(t)~\doteq~\max\{ k\geq 0\,;~~t_k^\ve\leq t\},$$
in view of  (\ref{RR}) we can write
\bel{final}\bega{l}\ds
\bfk(t,s)-\bfk(t_0,s)
~=~
\lim_{\ve\to 0}
\sum_{k=0}^{k^\ve(t)}
 \left(\ve\int_0^{t_k\wedge s} \Psi\bigl(t^\ve_k,\sigma, P^\ve(t_k^\ve,\sigma), 
 \bfk^\ve(t_k^\ve,\sigma)\bigr)d\sigma
\right)\times\bfk^\ve(t_k^\ve,s)
\\[4mm]
\ds
+\lim_{\ve\to 0}\sum_{k=0}^{k^\ve(t)}
\Big(\int^{t_k}_{s}e^{-\beta(t^\ve_k-s)}\nabla\Phi(P^\ve(t^\ve_k,s'))\times
(P^\ve(t^\ve_{k-1},s')-P(t^\ve_{k-1},s))d\mu^\ve_{k}(s')\Big)
\times 
\bfk^\ve(t_k,s).\enda\eeq
By the uniform convergence (\ref{conv}) and the smoothness of the function $\Phi$,
the first term on the right hand side of (\ref{final}) 
converges to the corresponding term in (\ref{E6}).
Moreover, recalling (\ref{medef}), 
and using the weak convergence $\mu^\ve\rightharpoonup\mu$, we conclude that 
the second term on the right hand side of (\ref{final}) converges
to the corresponding term in (\ref{E6}).
This completes the proof.
\endproof
\v
\section{More general models}
\label{s:6}
\setcounter{equation}{0}

All previous results can be extended to the case where 
$0<\alpha<+\infty$, so that all sections of the stem undergo
a linear elongation, exponentially decreasing in time.  
In addition, following  Remark~1,
one can also consider deformations which minimize the more general 
deformation energy (\ref{E33}), where the twist and bending components
are given different weights. 
We describe below the minor differences
in the analysis, required by these extensions.

When $\alpha< +\infty$, the unit tangent vector $\bfk$ to the stem is 
given by (\ref{k}), and 
the formula (\ref{P}) is replaced by (\ref{PP}).
The evolution of  $\bfk(t,s)$ is still described by 
 (\ref{E6}), replacing  (\ref{PSD}) with 
\bel{PSD2}
\Psi(t,\sigma, P, \bfk)~\doteq~\bigl(1-e^{-\alpha(t-\sigma)}\bigr)e^{-\beta(t-\sigma)}\Big( \kappa(\bfk\times \bfe_3)
+ (\nabla\psi(P)\times \bfk)\Big).\eeq
\v
Next, we derive the appropriate replacements for the cones 
$\Gamma$  in (\ref{RC5}) and $\Gamma'$ in (\ref{RC7}).

As in (\ref{k}), call $\bfk(t,\sigma)$ the unit tangent vector to the 
curve $P(t,\cdot)$ at the point $\sigma$. 
Given any  vector function $\omega(\cdot)$, we define
the orthogonal decomposition
$$\omega(\sigma)~=~\omega^{twist}(\sigma)+\omega^{bend}(\sigma),$$
where
$$\omega^{twist}(\sigma)~\doteq~\Pi_{\bfk(t,\sigma)}\omega(\sigma),
\qquad\qquad \omega^{bend}(\sigma)~\doteq~\Pi_{\bfk^\perp(t,\sigma)}\omega(\sigma),$$
denote the components parallel and orthogonal to the vector $\bfk(t,\sigma)$,
respectively.

 Consider the constrained optimization problem
\bel{E35}\hbox{minimize:}\quad \E~\doteq~{ 1\over 2}\int_0^{s'} 
\bigl(1-e^{-\alpha(t-\sigma)}\bigr)
e^{\beta(t-\sigma)}
 \Big(c_1\,\bigl|\omega^{twist}(\sigma)\bigr|^2+c_2
 \bigl|\omega^{bend}(\sigma)\bigr|^2
 \Big)\,d\sigma\,\eeq 
subject to the constraint
\bel{paw22}\nabla\Phi(P(t,s'))\cdot\left( \int_0^{s'}\bigl(1-e^{-\alpha(t-\sigma)}\bigr) \omega(\sigma)\times 
\bigl(P(t,s')-P(t,\sigma)\bigr)\, d\sigma\right)+\Phi(P(t,s'))~=~0. \eeq
To derive the appropriate necessary conditions,
consider a family of perturbations of the form
$$ \omega_{\epsilon}(\sigma)~=~
\omega(\sigma)+\epsilon \tilde{\omega}(\sigma)~=~\Big(
\omega^{twist}(\sigma)+\omega^{bend}(\sigma)
\Big)+\Big(\epsilon \tilde
\omega^{twist}(\sigma)+\epsilon \tilde\omega^{bend}(\sigma)
\Big).$$
Arguing as in (\ref{nc1}), 
we differentiate $\E$ w.r.t.~$\epsilon$ at $\epsilon=0$, 
and eventually obtain
\bel{nec_cnd1}\bega{l}\ds
\int_{0}^{s'}(1-e^{-\alpha(t-\sigma)})e^{\beta(t-\sigma)}\Big(c_{1}\omega^{twist}(\sigma)\cdot\tilde{\omega}^{twist}(\sigma)+  c_{2}\omega^{bend}(\sigma)\cdot \tilde{\omega}^{bend}(\sigma)\Big)d\sigma\\[4mm]
\ds\qquad 
=~\lambda \int_{0}^{s'}(1-e^{-\alpha(t-\sigma)})\tilde\omega(\sigma)\cdot\bigl(\nabla\Phi(P(t,s'))\times(P(t,s')-P(t,\sigma)) \bigr)d\sigma
\enda\eeq
for some Lagrange multiplier $\lambda$.
Notice that, by orthogonality, the inner products satisfy
$$\omega^{twist}(\sigma)\cdot\tilde{\omega}^{twist}(\sigma)
~=~\omega^{twist}(\sigma)\cdot \tilde\omega(\sigma)\qquad\qquad
\omega^{bend}(\sigma)\cdot \tilde{\omega}^{bend}(\sigma)~=~\omega^{bend}(\sigma)\cdot \tilde{\omega}(\sigma).$$
Since the relation (\ref{nec_cnd1}) holds true for 
every perturbation $\tilde{\omega}$,  one obtains
\bel{nec2}
c_{1}\omega^{twist}(\sigma)+c_{2}\omega^{bend}(\sigma)~=~
\lambda e^{-\beta(t-\sigma)}\nabla\Phi(P(t,s'))\times(P(t,s')-P(t,\sigma)).
\eeq
Projecting the above equation on the subspaces 
parallel and perpendicular to $\bfk(t,\sigma)$ respectively, with obvious meaning of notation we
finally obtain
\bel{nec3}\left\{\bega{rl}\ds
\omega^{twist}(\sigma)&\ds=~\frac{\lambda}{c_{1}}e^{-\beta(t-\sigma)}
\Big(\nabla\Phi(P(t,s'))\times(P(t,s')-P(t,\sigma)) \Big)^{twist},\\[4mm]
\ds 
\omega^{bend}(\sigma)&\ds=~\frac{\lambda}{c_{2}}e^{-\beta(t-\sigma)}
\Big(\nabla\Phi(P(t,s'))\times(P(t,s')-P(t,\sigma)) \Big)^{bend}.
\enda\right.\eeq
In place of (\ref{RC}), 
the {\bf cone of admissible velocities} can now be defined as
\bel{RC_mod} \bega{l}
\Gamma(t)~\doteq~\ds \Bigg\{ \bfv:[0,t]\mapsto \R^3\,;
~~\hbox{there exists a  positive measure $\mu$ supported on $\chi(t)$
such that}\\[4mm]
\quad \ds  \bfv(s)=-\int_0^s (1-e^{-\alpha(t-\sigma)})
e^{-\beta(t-\sigma)} \int_\sigma^t  
\bigg[\Big(\frac{\bfn(t,s'))}{c_{1}}\times(P(t,s')-P(t,\sigma))\Big)^{twist}
\\[4mm]\ds
\qquad \qquad\qquad \qquad
+\Big(\frac{\bfn(t,s')}{c_{2}}\times(P(t,s')-P(t,\sigma))\Big)^{bend} 
\bigg]   d\mu(s')  \times \bigl(P(t,s)-P(t,\sigma)\bigr)\, d\sigma \Bigg\}.
\enda
\eeq
Accordingly, the cone $\Gamma'(\tau)$ in (\ref{RC7}) is now replaced by
\bel{RC9}\bega{l}
\Gamma'(\tau)~\doteq~\ds \Bigg\{ \bfh:[0,T]\mapsto \R^3\,;~~
\hbox{there exists a  positive measure $\mu$, supported on
$\chi(\tau)$, such that} \\[4mm]
\ds \qquad
 \bfh(s)~=~ -\int_0^s (1-e^{-\alpha(t-\sigma)})e^{-\beta(t-\sigma)} \bigg(\int_\sigma^t  
\bigg[\Big(\frac{\bfn(t,s'))}{c_{1}}\times(P(t,s')-P(t,\sigma))\Big)^{twist}
\\[4mm]\ds
\qquad \qquad\qquad \qquad
+\Big(\frac{\bfn(t,s')}{c_{2}}\times(P(t,s')-P(t,\sigma))\Big)^{bend} 
\bigg]   d\mu(s') \bigg)d\sigma \times \bfk(t,s)\, d\sigma \Bigg\}.
\enda
\eeq

All the previous arguments can be applied to this  more general situation, 
with minor changes.

\v
\section{Numerical simulations}
\label{s:7}
\setcounter{equation}{0}

We present a couple of numerical simulations for the model 
\eqref{2d} in two dimensional space. 
Finite difference discretization is used, where the ODE \eqref{2d}
is solved with forward Euler method, and the integrations are
carried out with trapezoid rule. The stem/vine is discretized 
with uniform arc-length $\Delta s$, and the time step is also uniform 
with $\Delta t=\Delta s$. 
Simulations are carried out in Matlab. 
All the Matlab codes used,  together with many figures and simulations
can be found in \cite{WSweb}. 

\begin{figure}[ht]
%\centerline{\hbox{\includegraphics[clip,trim=13.7mm 10mm 10mm 9mm,width=6cm]{FIG/stemFigB.png}}
%\hbox{\includegraphics[clip,trim=13mm 10.5mm 10mm 9mm,width=6cm]{FIG/stemFigA.png}}}
\centerline{\hbox{\includegraphics[width=7cm]{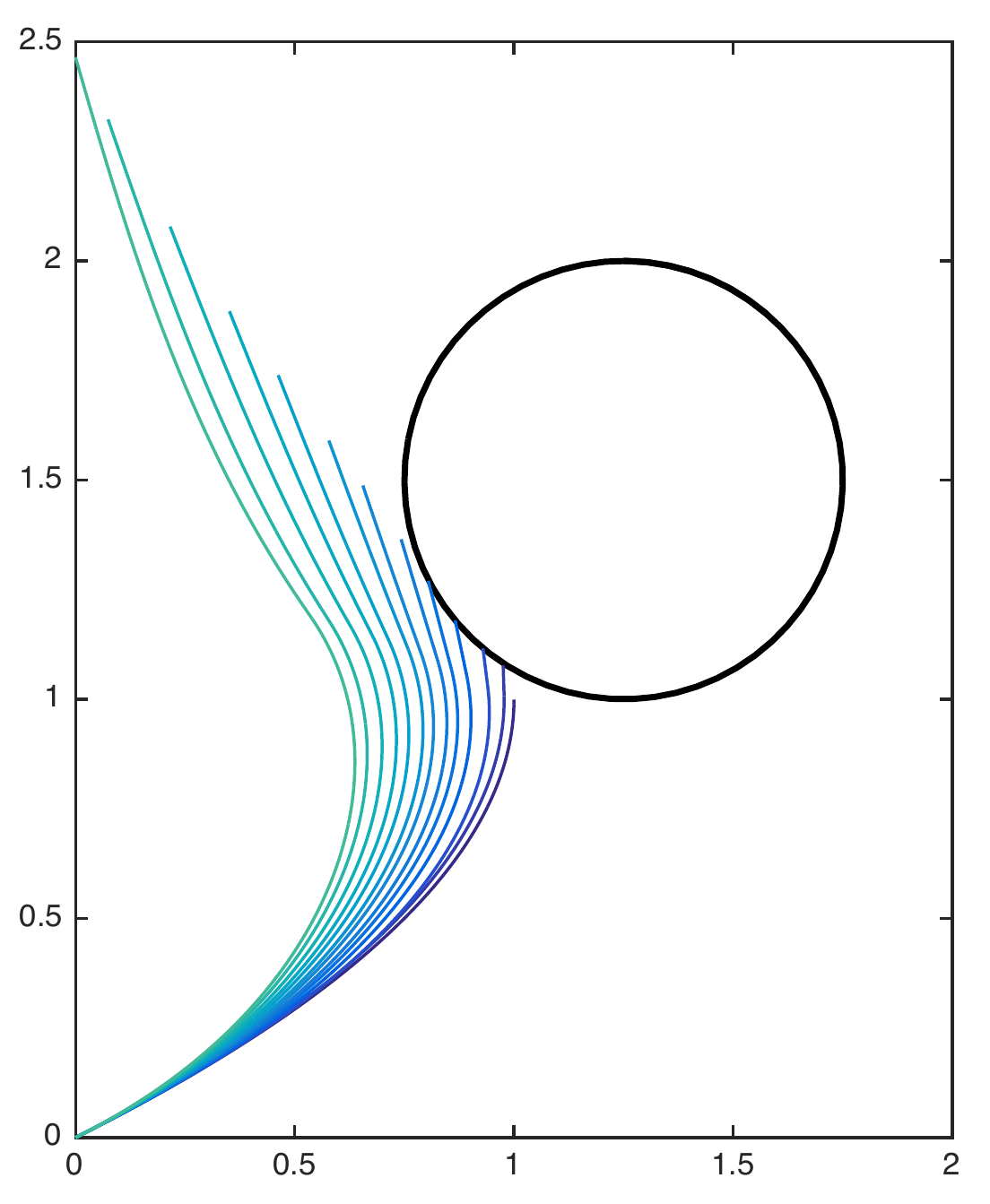}}
\hbox{\includegraphics[width=7cm]{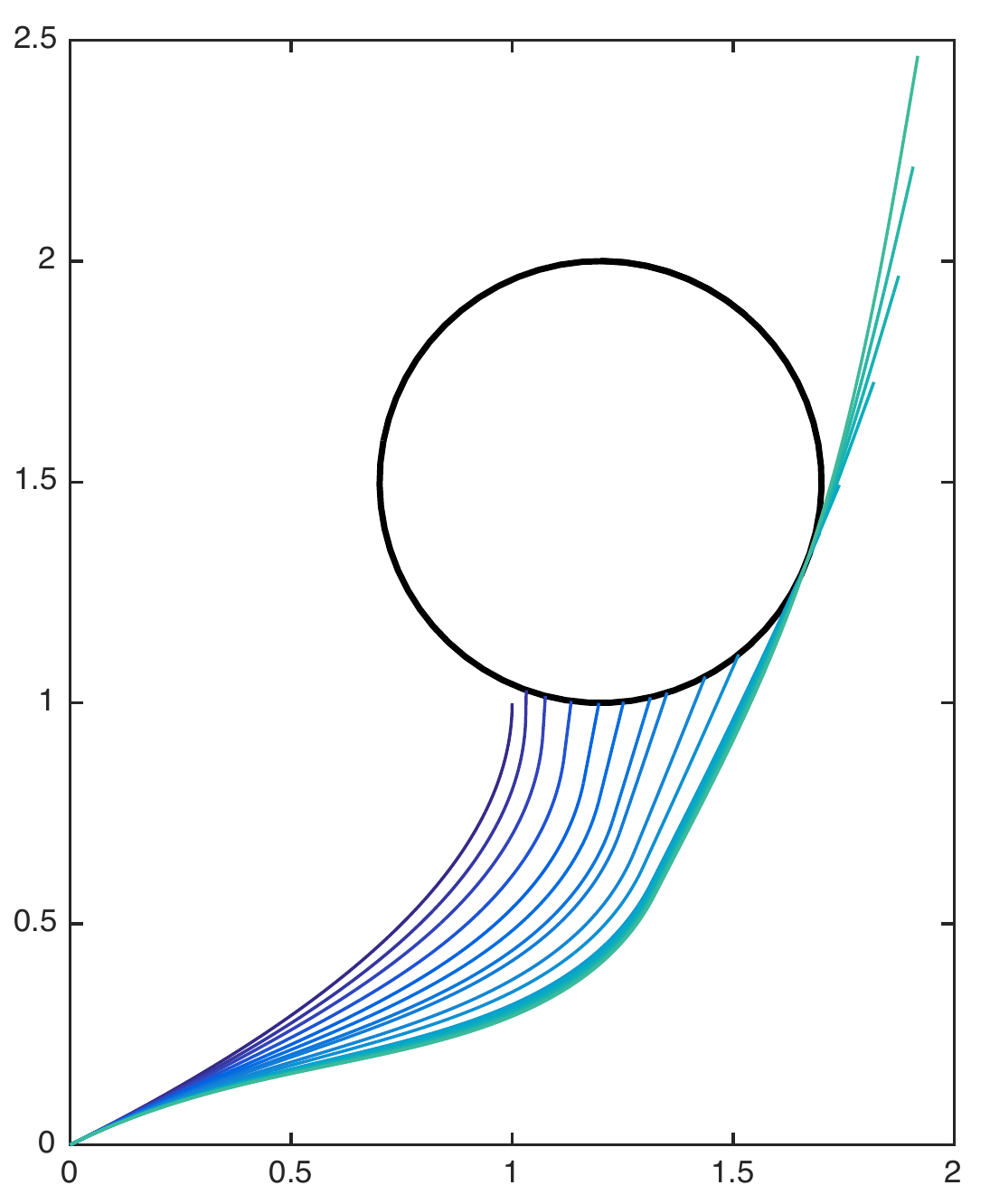}}}
\caption{\small Numerical simulation for testing the bifurcation for avoiding obstacles.}
\label{f:stemAB}
\end{figure}

\textbf{Simulation 1.} We simulate the growth of a tree stem avoiding an
obstacle along the way. 
In this model we neglect the second term on the right 
hand side of \eqref{2d}.
The following parameters are used:
$$ \beta=0.5, \quad \kappa=1.$$ % , \quad \Delta s=0.005.$$ 
The stem is initiated at the origin, with the initial shape 
$x=1-(y-1)^2$ for $0\le y\le 1$.  The obstacle is a circle, centered at 
$(a,b)$ with radius $r=0.5$. 
Two slightly different locations of $(a,b)$ are chosen, and the
results are shown in Figure~\ref{f:stemAB}. 
For the left plot, we use $(a,b)=(1.2,1.5)$, and the stem bends to the left
to avoid the obstacle. 
For the right plot, we use $(a,b)=(1.25,1.5)$, and the stem bends to the 
right to avoid the obstacle.

\begin{figure}[ht]
\centerline{\hbox{\includegraphics[clip,trim=8mm 6mm 4mm 5mm,height=8cm]{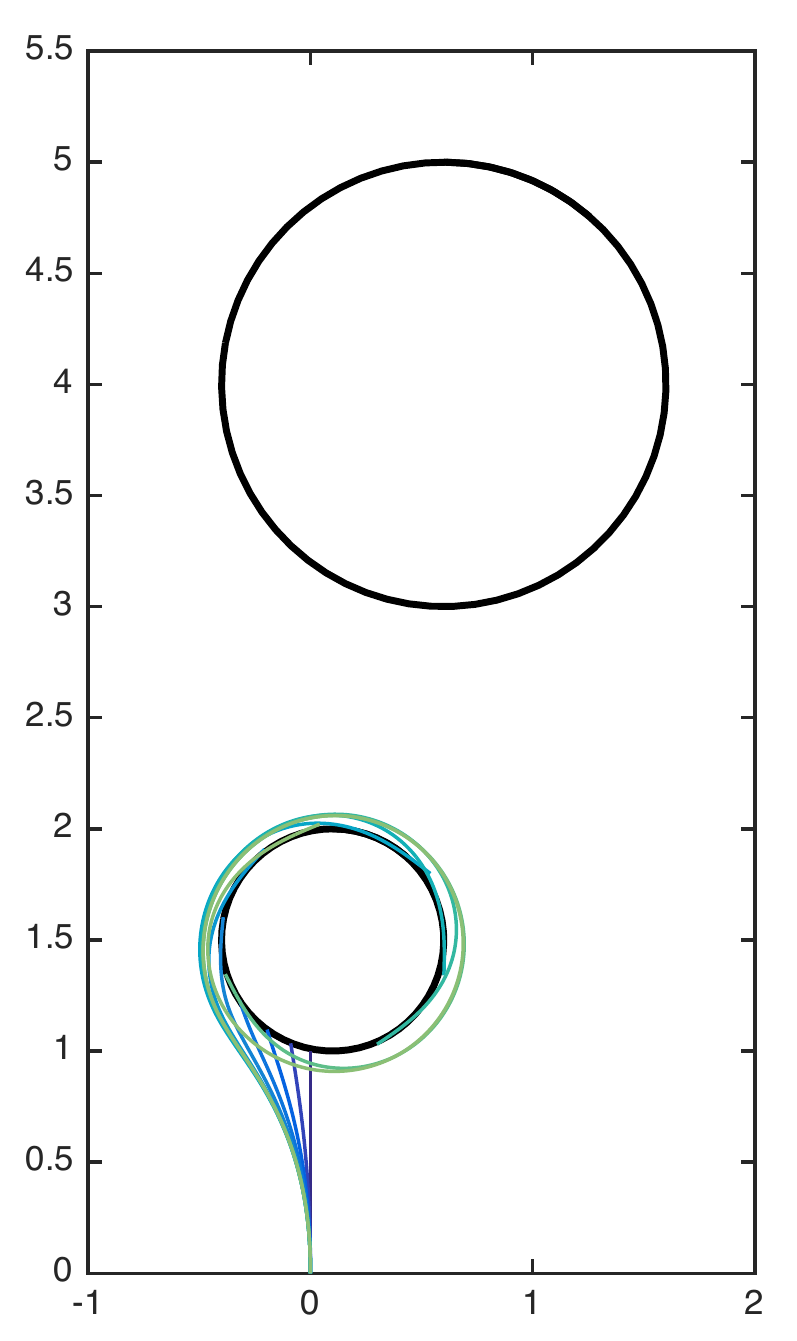}}
\quad\hbox{\includegraphics[clip,trim=8mm 5mm 4mm 4mm,height=8cm]{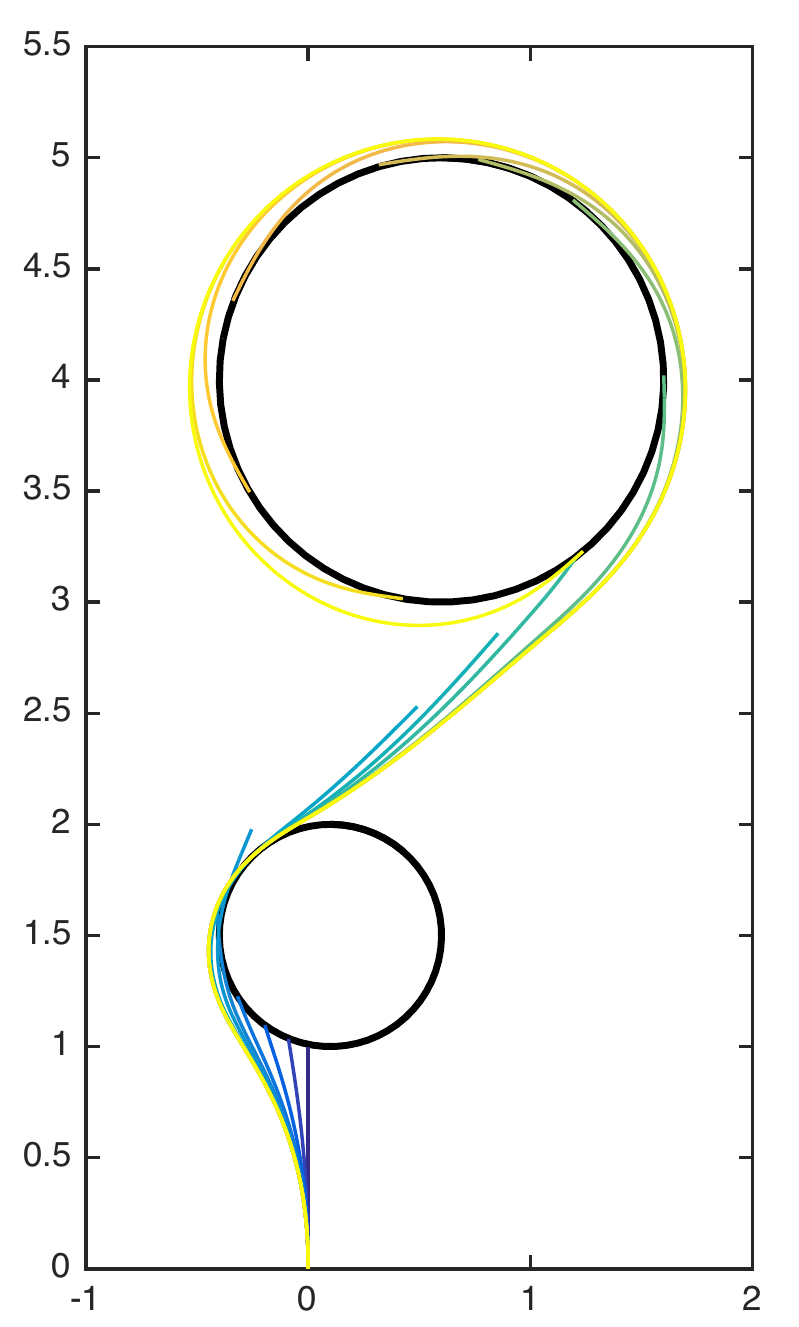}}
\quad\hbox{\includegraphics[clip,trim=8mm 7mm 4.5mm 4mm,height=8cm]{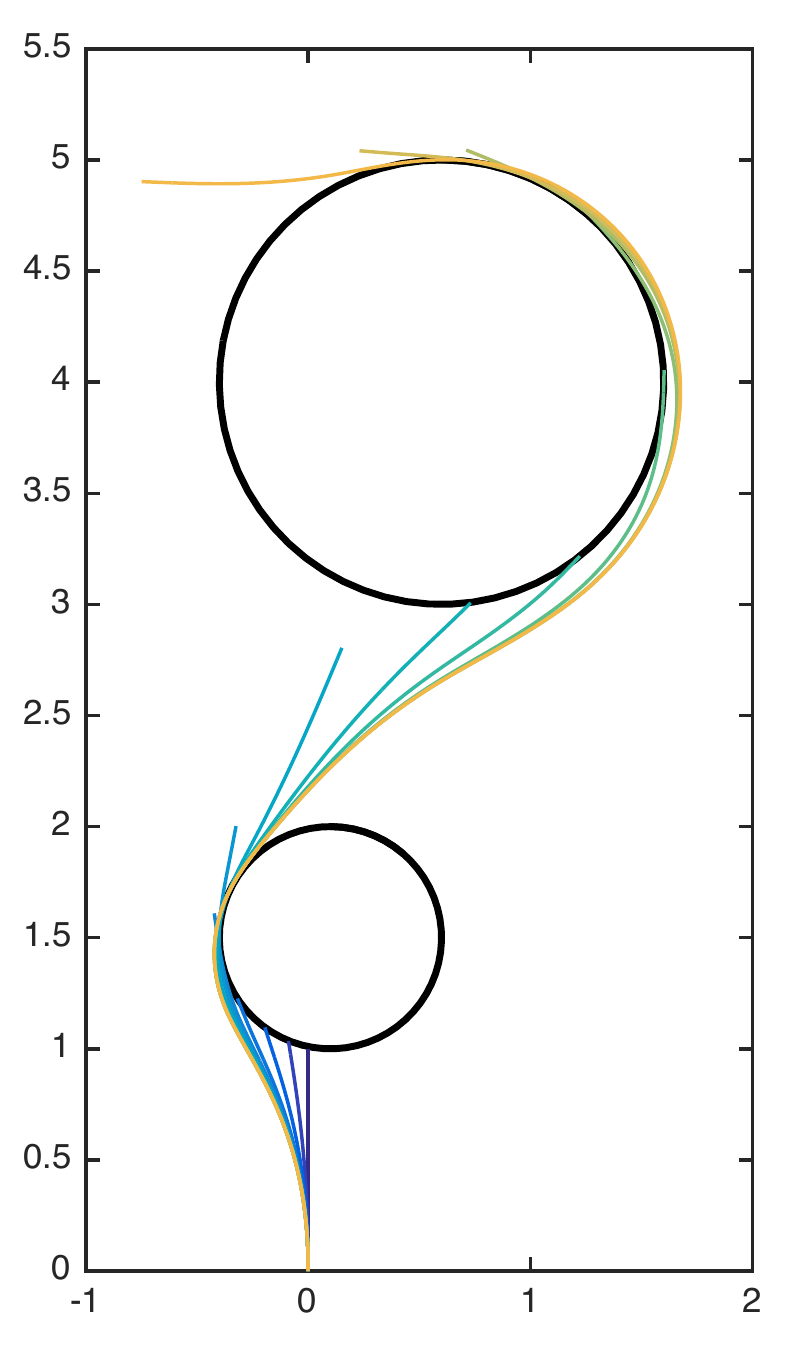}}}
\caption{\small Numerical simulation for vine growth with different bending 
parameter.}
\label{f:vineABC}
\end{figure}

\textbf{Simulation 2.}
We simulate a growing vine, in the presence of two obstacles.
The full equation \eqref{2d} is now used. 

The function $\eta(d)$ in \eqref{etaprop} is chosen to be 
$$ \eta(d) = \begin{cases} \gamma \left(1-e^{-d}\right), & 0\le d < \delta_0,\\
\gamma \left(1-e^{-\delta_0}\right), & d \ge \delta_0.
\end{cases}$$
The parameter $\gamma$ measures the strength of the 
feedback response, in the presence of an obstacle. 

We simulate various cases, highlighting the differences in the
solution caused by this bending factor $\gamma$. 
The vine is initiated at the origin, initially growing straight up.
Two circular obstacles are placed, one centered at $(0.1,1.5)$
with smaller radius $r_1=0.5$, the other centered at $(0.6,4)$ with 
larger radius $r_2=1$. 
We use the following parameters: 
$$ \beta=2, \qquad \kappa=1, \qquad \delta_0=0.05.$$

Numerical results with three different values of $\gamma$
are plotted in Figure~\ref{f:vineABC}:
\begin{itemize}
\item For the left plot, we use $\gamma=7$, a rather large value. 
The vine already curls around the smaller disc.
\item For the middle plot, we choose the smaller value 
$\gamma=4$.
The vine now fails to cling to the smaller disc, 
but manages to curl around the  larger disc.
\item For the right plot, we choose an even smaller parameter value: 
$\gamma=3$.
In this case,  the response to gravity prevails
and the vine fails to cling even to the larger disc.
\end{itemize}

\v

\end{document}